\def\|{||}
\def\qed{\hfill $\Box$}
\documentclass[11 pt]{amsart}
\usepackage{amssymb}
\usepackage{amsmath}
\usepackage{amsthm}

\theoremstyle{plain}
\newcommand{\cc}{{C_c^{\infty}}}
\newcommand{\mat}[4]{\left(\begin{array}{cc}#1&#2\\#3&#4\end{array}\right)}
\newcommand{\ppi}{{\tilde{\pi}}}
\newcommand{\G}{\tilde{G}}
\newcommand{\oo}{{\mathcal O}}
\newcommand{\ab}{{\underline{a}}}

\renewcommand{\L}{\tilde{L}}
\newcommand{\w}{\tilde{w}}
\newcommand{\Lv}{\L_{v}}

\renewcommand{\aa}{{\bf A}}
\newcommand{\W}{\tilde{W}}

\newcommand{\pp}{\tilde{\varphi}}

\renewcommand{\u}[1]{\underline{#1}}
\renewcommand{\v}[1]{\underline{\underline{#1}}}
\newtheorem{theorem}{Theorem}[section]
\newtheorem{proposition}[theorem]{Proposition}
\newtheorem{lemma}[theorem]{Lemma}
\newtheorem{corollary}[theorem]{Corollary}

\newtheorem{conjecture}[theorem]{Conjecture}
\theoremstyle{definition}
\newtheorem{remark}[theorem]{Remark}

\numberwithin{equation}{section}

\begin{document}
\pagenumbering{arabic}
\title[Central value]{Central value of automorphic
$L-$functions}
\author{Ehud Moshe Baruch}
\address{Department of Mathematics\\
         University of California, Santa Cruz\\
         Santa Cruz, CA 95064}
\email{baruch@math.ucsc.edu}
\author{Zhengyu Mao}
\address{Department of Mathematics and Computer Science\\
         Rutgers university\\
         Newark, NJ 07102-1811}
\email{zmao@andromeda.rutgers.edu}

\subjclass{Primary: subject; Secondary: subject}
\date{}

\begin{abstract}
We prove a generalization to the totally real field case of
the Waldspurger's formula relating the Fourier
coefficient of a half integral weight form and the central value of the
L-function of an integral weight form. Our proof is based on a new
interpretation of Waldspurger's formula in terms of equality
between global distributions. As applications we generalize
the Kohnen-Zagier formula for holomorphic forms and prove the equivalence
of the Ramanujan conjecture for half integral weight forms and a case of the
Lindelof hypothesis for integral weight forms. We also study the Kohnen
space in the adelic setting.
\end{abstract}
\maketitle
\section{Introduction}

In this paper, we generalize Waldspurger's formula 
 to the totally real field case, and study some of its applications. 
We use here the term {\em Waldspurger's formula} to mean an identity relating the Fourier
coefficients of a half integral weight form and the  central twisted $L-$values of an 
integral weight form. A key point in this paper is a new interpretation of Waldspurger's
formula. We study the formula in the adelic setting. This setting allows us to interpret
Waldspurger's formula as an equality between two global (adelic)
distributions.
Roughly, we define two global  distributions $I$ and $J$, and obtain the following
factorization
into products as  distributions over local fields:
$$I=c_1 \prod_v I_v,\,\,\,J=c_2\prod_v J_v.$$
Here $c_1,c_2$ are two global contants which can be interpreted respectively as
the central $L$-value of an integral weight form and as the square of 
the Fourier coefficient of a half integral weight form. 
For the more precise formulas, see
Propostions \ref{P:41} and \ref{P:42}.
Waldspurger's formula is then just the identity
$c_1=c_2$, (or more precisely, a family of such identities). Interpretted this way, 
Waldspurger's formula follows immediately from the comparison between the
global distributions $I$ and $J$ and the comparison
 between the local distributions $I_v$ and $J_v$. The formula
fits into a more general family of formulas that should result
from the comparison of the global distributions. Works on the theory of relative trace formula
have proved or conjectured many such comparisons of global distributions, 
see \cite{gu}, \cite{j2}, \cite{m}, \cite{mr} etc.
The resulting identities from such comparisons in
all these cases should also have arithmetic interest as in the case of Waldspurger's formula.

We look at some applications of the 
generalization of Waldspurgers formula, among them the equivalence of the Ramanujuan conjecture
for half-integral weight forms and a special case of
the Lindel\H{o}f hypothesis, (this result is useful in the work of Cogdell-Piatetski-Shapiro-Sarnak
on ternary forms \cite{cpss}). Our result is stated in the 
general setting of automorphic forms. As many of its applications are in terms of
holomorphic modular forms over ${\Bbb Q}$, we also translate our result into this
language. In the process of the translation, we generalize the well known 
Kohnen-Zagier formula, with the restrictions
on the quadratic twist dropped (Theorem~\ref{T:KZgeneral}). 

With our interpretation of the Waldspurger's formula,
 it is clear how the proof of the identity will proceed. The equality
of the global distributions $I$ and $J$ follows from the global theory of relative trace formula 
(\cite{j}), 
while the equalities of the local distributions follow from the corresponding
local theory, \cite{bm},\cite{bm2}.
The subtle part is to show that we can fit the identities $c_1=c_2$ into a
family of identities. There we need to use the beautiful results of Waldspurger on  theta
correspondence.

Below we describe the content of the paper in more detail.

\subsection{An explicit version of Waldspurger's formula}

The Shimura correspondence associates
a cusp form $f(z)$ with integral
weight $2k$ to a half integral weight cusp
form $g(z)$ with weight $k+1/2$. 
Waldspurger is the first to described a relation
between the twisted central $L-$values $L(f,D,k)$ and the Fourier coefficients
of $g(z)$,
 \cite{wa1}. There are  many later versions of Waldspurger's formula.
 See for example the papers of Shimura \cite{s},\cite{sh},
 Niwa
\cite{n}, Katok-Sarnak \cite{ks}, Kohnen-Zagier \cite{k},\cite{k1},
Gross \cite{gr}. 
  The result of  Kohnen and Zagier is probably the easiest to describe, it states: \cite{k1}.

\begin{quote}
Let $f(z)$ be a new form of weight $2k$, square free and odd
level $N$, of trivial character.
There is a unique (up to a scalar multiple) weight $k+1/2$ form $g(z)=\sum
_{n=1}^{\infty}c(n)e^{2\pi inz}$  corresponding to $f$
and lying in the {\em Kohnen Space} (\cite{k2}), such that when
 $D$ is a fundamental discriminant with
 $(-1)^k D>0$ and $(\frac{D}{l})=w_l$  for all prime divisors $l$ of
$N$, 
\begin{equation}\label{E:11}
\frac{|c(|D|)|^2}{<g,g>}=\frac{(k-1)!}{\pi^k}|D|^{k-1/2}2^{\nu(N)}
\frac{L(f,D,k)}{<f,f>}
\end{equation}
\end{quote}

In the above statement, we have adopted the notations in \cite{k1}, and
$\nu(N)$ is the number of prime divisors of $N$, $w_l$ is
the eigenvalue of the Atkin-Lehner involution acting on $f(z)$.

For the many applications of an equation of type (\ref{E:11}), see
for examples \cite{iw1},\cite{k},\cite{os},\cite{lr}.

\subsection{Waldspurger's formula: adelic version}

In this paper, we will work with the more general notaion of automorphic representations
over a totally real number field $F$.
The relationship between the modular forms and the automorphic representations is
as follows.
 A half-integral 
weight modular form is a vector in the space of
an automorphic representation $\tilde{\pi}$
on $\widetilde{SL}_2$, the double cover of $SL_2$. An integral weight 
modular form with trivial character can be considered as a vector in the space
of an automorphic representation $\pi$ on $PGL_2$. 
The representation theoretic version of the Shimura
correspondence is a theta correspondence relating $\pi$ on $PGL_2$
to some $\tilde{\pi}$ on $\widetilde{SL}_2$, \cite{wa}.

 Our first task
is to define the constants associated to $\tilde{\pi}$ and $\pi$
that are the analogues of the Fourier coefficients of the modular forms.
Such constants are defined in \S~\ref{S:constants}. For $D\not=0\in F^*$,  we define
the $D-$th ``Fourier coefficients" of $\pi$ and $\ppi$ to be
$d_{\pi}(S,\psi^D)$ and $d_{\ppi}(S,\psi^D)$, (see \S~\ref{S:constants} for
notations). Our version of the Waldspurger's formula is  (Theorem~\ref{T:21}): for
a $\pi$ and $D$, there is a corresponding representation $\ppi=\ppi(\pi,D)$  of 
$\widetilde{SL}_2$, such that
\begin{equation}\label{E:Waldnew}
|d_{\pi}(S,\psi^D)|^2 L^S(\pi,1/2)=|d_{\ppi}(S,\psi^D)|^2.
\end{equation}
Here $L^S(\pi,1/2)$ is the central (partial) $L-$value.

\subsection{The role of Waldspurger's result on theta correspondence}

A clear difference between our formula and (\ref{E:11}) is that the twisted
$L-$value does not appear in our formula. To introduce the twisted $L-$value,
we apply the formula to $\pi\otimes \chi_D$ where $\chi_D$ is a quadratic
character associated to $D\in F^*$. 
The problem of course is that $\ppi(\pi\otimes\chi_D, D)$ is not
necessarily the same  for all $D$.
 
 The results of Waldspurger on theta correspondence (\cite{wa2})
 gives a partition of the set 
$F^*$ into a finite collection of subsets, such that  the representation
 $\ppi(\pi\otimes\chi_D, D)$ remains the same  for $D$ in a given subset.
 The equation (\ref{E:Waldnew})  then gives
a formula for the twisted $L-$value $L(\pi\otimes \chi_D,1/2)$ in terms
of the $D-$th Fourier coefficients of $\ppi$, as long as 
$D$ lies in this particular subset. 

\subsection{Generalization of the Kohnen-Zagier formula}

We compare our result with the Kohnen-Zagier formula
(\ref{E:11}). We now understand the conditions on $D$ in 
the Kohnen-Zagier formula. The condition is to
ensure that $D$ lies in a given subset of ${\Bbb Q}^*$, so
that the  half-integral weight form appear in (\ref{E:11}) 
remains the same. With this understanding, we see that for $D$
in other subsets of ${\Bbb Q}^*$, there should also be another
version of the Kohnen-Zagier formula, involving a different
half integral weight form.

In our generalization, we assoicate to $f(z)$ a partition of 
${\Bbb Q}^*$, and to each subset $X$ of the partition a 
half integral weight forms $g_X(z)$. We get an equation in the form
of  (\ref{E:11}) for each $g_X(z)$, which holds for all fundamental discriminants $D$
that lie in the subset $X$.

In the process of deriving the generalized Kohnen-Zagier
formula from our adelic version, we give the adelic
interpretation of the concept of Kohnen space. To a half-integral weight
eigenform of level $4N$ ($N$ odd)
associates a vector $\pp$ in the space of an automorphic representation
$\ppi=\otimes\ppi_v$ of $\widetilde{SL}_2$, where $\pp=\pp_{\infty}\otimes\pp_2\otimes\pp_3\otimes\ldots$ 
with $\pp_v$ being vectors in the space $V_{\ppi,v}$ of $\ppi_v$. The
vector $\pp_2$
could lie in a two dimensional subspace of $V_{\ppi_2}$. The Kohnen
space is a subspace of half-integral weight forms. We show that it is exactly the subspace generated
by $\pp$'s whose local component
$\pp_2$ lies in a particular one dimensional subspace of $V_{\pi,2}$.

{\bf Structure of the paper}:

The paper is organized as follows: In \S~\ref{S:constants}, we define the constants $d_{\pi}(S,\psi)$
and $d_{\ppi}(S,\psi)$. We recall Waldspurger's result on theta correspondence
in \S~\ref{S:Waldspurger}. In \S~\ref{S:mainresult} we state the main theorems.
The relative trace formula of \cite{j} is reviewed in \S~\ref{S:Jacquet}. We describe
the local theory of the relative trace formula in \S~\ref{S:BM}. 
The proof of the
main theorems are  given in \S~\ref{S:proof}.
In \S~\ref{S:local},  we compute some examples of local constants
appearing in the identity for $L(\pi,1/2)$. The computations are
just easy exercises, and the results are used in the translation from
adelic language to modular form language of our formula, as well as in a proof
in \S~\ref{S:proof}.  In \S~\ref{S:dictionary}, we give a 
dictionary between the language of representation theory and modular forms.
We also give a discussion on  the Kohnen space.
In \S~\ref{S:KZ}, we
prove the Kohnen-Zagier formula without the conditions on
the fundamental discriminant $D$. 

{\bf Notations and background}:

Let $F$ be a totally real
number field, $\aa$ its adele ring. We will use $v$ to
denote places of $F$. When $v$ is nonarchimedean, let $\oo_v$ be the ring of integers in $F_v$,
$P_v$ (or $P$) be its prime ideal, $\varpi$ its uniformizer,
and $q_v$ (or $q$) the size of the residue field $\oo_v/P$.
We use $||_v$ to denote the metric on $F_v$.

Let $G=GL_2$,
$\G=\widetilde{GL}_2$ and $G'=\widetilde{SL}_2$. We will use $e$ to denote
the identity elements of the groups $G$, $G'$ and $\G$.
Let $Z$ be the center of $G$. Then $PGL_2=G/Z$.
Let
$B$ be the subgroup of $GL_2$ consisting of upper triangular
matrices, $\tilde{B}$ be its lifting to $\G$.

We will use $(g,\pm 1)$ to denote an element in $\G$. Let $[*,*]$ be the Hilbert symbol.
The multiplication in $\G$ takes the form:
$$(g_1,1)\cdot(g_2,1)=(g_1g_2,[\frac{x(g_1g_2)}{x(g_1)},
\frac{x(g_1g_2)}{x(g_2)}\det g_1])$$
where for $g=\mat{a}{b}{c}{d}$,
$x(g)=c$ if $c\not=0$ or $d$ if $c=0$.

Let $n(x)=\mat{1}{x}{}{1},\tilde{n}(x)=(n(x),1).$ Let
$w=\mat{}{1}{1}{}$, $\tilde{w}=(\mat{}{-1}{1}{},1).$  Let $\u{a}=
\mat{a}{}{}{1}$, $\v{a}=(\mat{a}{}{}{a^{-1}},1)$.

We fix a nontrivial  additive character $\psi$ of $\aa/F$. Then $\psi=\prod_v\psi_v$.
 For $D\in F^{*}$,
let $\psi^D(x)=\psi(Dx)$; let $\chi_D$ be the quadratic character
of $\aa^*/F^*$ associated to the field extension $F(\sqrt{D})$.
At a local place $v$, for $D\in F^*_v$, we let $\chi_D$
be the quadratic character of $F^*_v$ associated to the extension
$F_v(\sqrt{D})$.

We will fix measures as follows.  The choice of additive $dx$ on $F_v$
measure does not matter for the statement of our theorem. We will however
fix it
to be self dual for the character $\psi_v$.
 The multiplicative measure is $d^{*}a=
(1-q^{-1})^{-1}\frac{da}{|a|_v}$, where $q$ is the size of the residue field when $v$
is $p-$adic, and $q=\infty$ when $v=\infty$.
We fix the measures  for $Z\backslash GL_2$ and $\widetilde{SL}_2$ as in \cite{bm}
and \cite{bm2}. Write $g=z(c)n(x)w\u{a}n(y)$ for $g\in G(F_v)-B(F_v)$, then
$dg=|a|_vd^{*}cd^{*}adxdy$ is the measure on $G(F_v)$. The measure on $Z(F_v)$
is $dz(c)=d^{*}c$, and we use the resulting quotient measure on $Z(F_v)\backslash G(F_v)$.
For $g\in G'(F_v)-\tilde{B}(F_v)\cap G'(F_v)$ with $g=\tilde{n}(x)\tilde{w}\v{a}\tilde{n}(y)$,
we define $dg=|a|_v^2d^{*}adxdy$ to be the measure on $G'(F_v)$.

Define the Weil constant $\gamma(a,\psi_v^D)$ over $F_v$ to satisfy:
$$\int\hat{\Phi}(x)\psi_v^D(ax^2)dx=|a|_v^{-1/2}\gamma(a,
\psi_v^D)\int\Phi(x)\psi_v(-a^{-1}x^2)dx$$
where
$$\hat{\Phi}(x)=\int \Phi(y)\psi_v^D(-2xy)dy.$$
We let $\tilde{\gamma}(a,\psi_v^D)=\frac{\gamma(a,\psi_v^D)}{
\gamma(1,\psi_v^D)}[-1,a]$.

We  use $\pi$ to denote an
irreducible cuspidal
representation of $G(\aa)$ with trivial central character, and
use $\ppi$ to denote an
irreducible cuspidal
representation of $G'(\aa)$. $\pi$ can
be considered as a representation of $PGL_2(\aa)$. We have
$\pi=\otimes \pi_v$ and $\ppi=\otimes\ppi_v$ 
as the restricted tensor products of representations over local
fields. 
We will use $V_{\pi}, V_{\ppi}, V_{\pi,v}$ and $V_{\ppi,v}$ to 
denote the spaces where
the representations $\pi,\ppi,\pi_v$ and $\ppi_v$ act on respectively..

When $\mu$ is a character of $F_v^{*}$, we will let
$\pi(\mu,\mu^{-1})$ denote the principal series representation
of $G(F_v)$ induced from $\mu$. It acts by right translation
on the space of functions
$\phi$ on $G(F_v)$ that satisfies:
\begin{equation}\label{E:12}
\phi(n(x)\u{a}zg)=\mu(a)|a|_v^{1/2}\phi(g)
\end{equation}
We use $\ppi(\mu,\psi_v)$ to denote the principal series
representation of $G'(F_v)$ that acts on the
space of functions $\phi$ of $G'(F_v)$
satisfying:
\begin{equation}\label{E:13}
\phi(\tilde{n}(x)\cdot\v{a}\cdot g)=\mu(a)\tilde{\gamma}(a, \psi_v)|a|_v\phi(g)
\end{equation}
These representations are unramified if $\mu$ and $\psi_v$
are unramified. 

The $L-$function $L(\pi,s)$ is defined in \cite{jl}.
So is the factor $\epsilon(\pi,s)=\prod \epsilon(\pi_v,s,\psi_v)$. We
note that $\epsilon(\pi_v,1/2)=\epsilon(\pi_v,1/2,\psi_v)$ is
independent of the choice of $\psi_v$.

By the well known result on the Shimura-Waldspurger (theta)
 correspondence (\cite{rs}, \cite{s}, \cite{wa2}), 
 for the given $\psi^D$ there associates
 a unique irreducible cuspidal
representation $\ppi(\pi, D)=\Theta(\pi,\psi^D)$ of $G'(\aa)$ (\cite{wa}).
Here $\Theta$ denotes the theta correspondence.  Similarly given $\ppi$ on $G'$,
there is a unique irreducible cuspidal representation $\Theta(\ppi,\psi^D)$
on $PGL_2$ (\cite{wa}). We note that the space of $\Theta(\pi,\psi^D)$
and $\Theta(\ppi,\psi^D)$ could be zero dimensional. The theta correspondence
is also defined locally, which we again denote by $\Theta$.

We will use $S$  to denote a finite set of local places containing all $v$
which is archimedean or has even residue charactersitic.
For $v\not\in S$, the covering $G'(F_v)$ splits
over $SL_2(\oo_v)$. With this splitting, we consider
$SL_2(\oo_v)$ a subgroup of $G'(F_v)$. Explicitly the embedding of
$SL_2({\mathcal O}_v)$  in $G'$ is given
by $g\mapsto (g,\kappa(g))$,
where
$\kappa(\mat{a}{b}{c}{d})=[c,d]$ if $|c|_v<1$ and $c
\not=0$, $\kappa(g)$ equals 1 when $|c|_v=1$ or $c=0$.

We will use $\|\varphi\|$ to denote the norm of a vector
$\varphi$: if $(*,*)$ is a Hermitian form on a space $V$, for
$\varphi\in V$, let $\|\varphi\|=(\varphi,\varphi)^{1/2}$.
 We use $\{\delta_i\}$ to denote a set of
representatives for the square classes in $F_v^{*}/(F^*_v)^2$,
with $\delta_1=1$.

{\bf Acknowledgement}: Professor Jacquet suggested to us to consider
the application of \cite{j} in $L-$functions. We thank him and J.
Cogdell, B.Conrey, S. Gelbart, D. Ginzburg, E. Lapid, S. Rallis, 
D. Ramakrishnan, P. Sarnak, D. Soudry,
J. Sturm 
 for helpful conversations. We would also like
to thank The Ohio State University, Institute for advanced study,
the Weizmann institute of science for their hospitality during the
visit of one or both authors, and Gelbart and Rallis in particular
for their invitation. The second author was partially supported
by NSF DMS 9304580.

\section{Definition of two constants}\label{S:constants}

\subsection{A constant associated to $\pi$ on $G$}\label{S:constant1}

We define a constant $d_{\pi}(S,\psi)$ which can be considered as the
Fourier coefficient of $\pi$. The constant depends only on the character
$\psi$, the choice of the finite set of places $S$ and the choice of
Haar measures. Note that the choice of Haar measures and $\psi$ are fixed
in the introduction.
 
\subsubsection{Whittaker model on G}

Let $\pi$ be an irreducible cuspidal automorphic representation
of $G(\aa)$ with trivial central character. 
 Let $V_{\pi}\subset L^2(Z(\aa)G(F)\backslash G(\aa))$ be the
space that $\pi$ acts on.
 For $\varphi\in V_{\pi}$, let
\begin{equation}\label{E:defW}
W_{\varphi}(g)=\int_{\aa/F}\varphi(n(u)g)\psi(-u)du.
\end{equation}
Then the space $\{W_{\varphi}|\varphi\in V_{\pi}\}$ gives
the global Whittaker model of $\pi$. 
\begin{remark}
When an integral weight form $f$ is considered as a vector $\varphi$
in the space
of $V_{\pi}$, its Fourier coefficients are roughly the values of $W_{\varphi}(e)$
for various choices of $\psi$, (see \S~\ref{S:dictionary}).
\end{remark}

  For any admissible
representation $\pi_v$ of  $G(F_v)$, the space of its
$\psi_v$-Whittaker functional $L_v: V_{\pi,v}\rightarrow {\Bbb C}$
satisfying:
\begin{equation}\label{E:21}
L_v(\pi_v(n(u))v)=\psi_v(u)L_v(v), v\in V_{\pi,v}
\end{equation}
is at most one dimensional. For the $\pi_v$'s appear as local components
of  $\pi$,
such a space is one dimensional.
We will fix the linear form $L_v$ for any given $\pi_v$.

Let $S$ be a finite set of places as in the introduction, and contain all places $v$ where
$\pi_v$ is not unramified.
For $v\not\in S$,  $\pi_v$ is an unramified representation
of $G(F_v)$; let $\varphi_{0,v}\in V_{\pi,v}$ be the unique vector
fixed under the action of $G(\oo_v)$ such that $L_v(\varphi_{0,v})=1$.

We note that
$$L:\varphi\rightarrow W_{\varphi}(e)$$
is a linear form on $V_{\pi}$ satisfying (\ref{E:21}). 
From the uniqueness
of the local Whittaker functionals,
$L$ can been expressed as a product of local linear forms $L_v$. There is a 
constant $c_1(\pi,S,\psi,\{L_v\})$,
such that whenever $\varphi=\otimes_{v\in S}\varphi_v\otimes_{v\not\in S}
\varphi_{0,v}$, (here we  fix an identification between $V_{\pi}$
and the restricted tensor product $\otimes V_{\pi,v}$ where $V_{\pi,v}$ is the space of the local component
$\pi_v$)
\begin{equation}\label{E:defc1}
W_{\varphi}(e)=c_1(\pi,S,\psi,\{L_v\})\prod_{v\in S}L_v(\varphi_v).
\end{equation}

\subsubsection{Hermitian forms on $G$}

Define a Hermitian form on $V_{\pi}$ by:
\begin{equation}\label{E:22}
(\varphi, \varphi')=\int_{Z(\aa)G(F)\backslash G(\aa)}
\varphi(g)\overline{\varphi'(g)}dg
\end{equation}
Over a local place $v$, for a unitary representation $\pi_v$
with a nontrivial Whittaker functional $L_v$, we can define
a $G_v$-invariant Hermitian form on $V_{\pi,v}$ by:
\begin{equation}\label{E:23}
(v,v')=\int_{F^{*}}L_v(\pi_v(\ab)v)\overline{L_v(\pi_v(\ab)v')}\frac{da}{|a|_v}
\end{equation}
(see [Go]). From the uniqueness of $G_v$-invariant
Hermitian forms, we get: there is a constant $c_2(\pi,S,\psi,\{L_v\})>0$,
such that whenever $\varphi=\otimes_{v\in S}\varphi_v\otimes_{v\not\in S}
\varphi_{0,v}$,
\begin{equation}\label{E:defc2}
\|\varphi\|=c_2(\pi,S,\psi,\{L_v\})\prod_{v\in S}\|\varphi_v\|.
\end{equation}

\subsubsection{The constant $d_{\pi}(S,\psi)$}

\begin{lemma}\label{L:defd1}
The constant $d_{\pi}(S,\psi)$ defined by
$$d_{\pi}(S,\psi)=|c_1(\pi,S,\psi,\{L_v\})/c_2(\pi,S,\psi,\{L_v\})|$$
is independent of the choice of the linear forms $L_v$.
\end{lemma}
\begin{proof}
From the uniqueness of local Whittaker functionals, any other
choice of linear forms $L'_v$ must have the form $L'_v=a_vL_v$ with
$a_v$ some nonzero complex constants. From the definition, we get
$$c_2(\pi,S,\psi,\{L'_v\})=\prod_{v\in S}|a_v|^{-1}c_2(\pi,S,\psi,\{L_v\}),$$
$$c_1(\pi,S,\psi,\{L'_v\})=\prod_{v\in S}a_v^{-1}c_2(\pi,S,\psi,\{L_v\}).$$
Thus the constant $d_{\pi}(S,\psi)$ is independent of the choice of $\{L_v\}$.
\end{proof}

This is the ``Fourier coefficient" we associate to $\pi$. The constant
$d_{\pi}(S,\psi)$ is well defined as we fixed the choice of $\psi$ and 
the measure on $G$, (it is easy to check that the constant is
independent of the choice of additive
measure). Explicitly for any vector $\varphi=
\otimes_{v\in S}\varphi_v\otimes_{v\not\in S}
\varphi_{0,v}$ with $L_v(\varphi_v)\not=0$ for $v\in S$,
\begin{equation}\label{E:d1explicit}
d_{\pi}(S,\psi)=\frac{|W_{\varphi}(e)|}{\|\varphi\|}\prod_{v\in S}
\frac{\|\varphi_v\|}{|L_v(\varphi_v)|}.
\end{equation}

We make an observation on the dependence of 
 $d_{\pi}(S,\psi)$ on $\psi$.
\begin{lemma}\label{L:dpi}
Let $D\in F^*$. If $S$ is large enough such that $|D|_v=1$ for all $v\not\in S$, then
$$d_{\pi}(S,\psi)=d_{\pi}(S,\psi^D).$$
\end{lemma}
\begin{proof}
Take a vector $\varphi=\otimes_{v\in S}\varphi_v\otimes_{v\not\in S}
\varphi_{0,v}$ in the space of $\pi$. 
We will let $L^D_v(\varphi_v)=L_v(\pi_v(\u{D})\varphi_v)$. Then
$L^D_v$ is a nontrivial $\psi_v^D$-Whittaker functional on $\pi_v$.
Let $\|\varphi_v\|_D$ be the norm of $\varphi_v$ defined by (\ref{E:23})
with $L_v$ replaced by $L^D_v$. 
When $v\not\in S$, clearly $L_v(\pi_v(\u{D})\varphi_{0,v})=L_v(\varphi_{0,v})=1$. Thus
using the above explicit form, 
$$d_{\pi}(S,\psi^D)=\frac{|W_{\varphi}^D(e)|}{\|\varphi\|}\prod_{v\in S}
\frac{\|\varphi_v\|_D}{|L^D_v(\varphi_v)|}.
$$
From a simple change of variable we get $W_{\varphi}^D(e)=W_{\pi(\u{D})\varphi}(e)$.
By the uniqueness of local Whittaker functional,
$$\frac{|W_{\varphi}^D(e)|}{\prod_{v\in S}|L^D_v(\varphi_v)|}
=\frac{|W_{\pi(\u{D})\varphi}(e)|}{\prod_{v\in S}|L_v(\pi(\u{D})\varphi_v)|}
=\frac{|W_{\varphi}(e)|}{\prod_{v\in S}|L_v(\varphi_v)|}.$$
 From
(\ref{E:23}), $\|\varphi_v\|_D=\|\varphi_v\|$. Thus we get the
equality in the Lemma.
\end{proof}

\subsection{A constant associated to $\ppi$ on $G'$}\label{S:constant2}

Let $\ppi$ be an irreducible cuspidal automorphic representation of $G'$.
We associate a constant $d_{\ppi}(S,\psi^D)$ to $\ppi$ in a similar
fashion. 
Let
$V_{\ppi}$ be the space of automorphic forms that $\ppi$ acts
on. 
For $\pp\in V_{\pi}$, let
$$\W_{\pp}^D(g)=\int_{\aa/F}\pp(\tilde{n}(x)\cdot g)
\psi^D(-x)dx.$$
Then  the Fourier coefficients of
half integral weight form can be intepreted as some $\W_{\pp}^D(e)$,
(see equation (\ref{E:fc2})).

We will assume $\ppi$ has a nontrivial $\psi^D$-Whittaker model, namely
$\W_{\pp}^D(g)\not=0$ for some $\pp\in V_{\ppi}$. Then locally, $\ppi_v$
has a nontrivial $\psi_v^D$-Whittaker model, unique up to a scalar multiple.
We will fix the corresponding linear forms $\L_v^D$ satisfying
for all $\pp_vv\in V_{\ppi,v}$,
$$\L_v^D(\ppi_v(\tilde{n}(x))\pp_v)=\L_v^D(\pp_v)\psi_v^D(x).$$

Let $S$ be a finite set of places
as in the introduction, and contain all places $v$ where
$\ppi_v$ is not unramified.
When $v\not\in S$, $\ppi_v$ is unramified and possess a nontrivial
$\psi_v^D$-Whittaker model $\L_v^D$; there is a unique vector in
$\pp_{0,v}$ that is fixed under $SL_2(\oo_v)$ and satisfying
$\L_v^D(\pp_{0,v})=1$.

The space $V_{\ppi}$ has the Hermitian form
\begin{equation}\label{E:25}
(\pp,\pp')=\int_{SL_2(F)\backslash G'(\aa)}\pp(g)\overline{\pp'(g)}dg.
\end{equation}
Over the space $V_{\ppi,v}$, one can define a Hermitian form
similar to (\ref{E:23}), though the definition is more complicated.
Let $\{\delta_i\}$ be a set of representatives of $F_v^{*}/
(F_v^{*})^{ 2}$, with $\delta_1=1$. From
\cite{bm},\cite{bm2}, we see there is a choice of $\psi_v^{D\delta_i}$-Whittaker
models $\L^{D\delta_i}_v$ (could be trivial) on $V_{\ppi,v}$, such that
$\L^{D\delta_1}_v=\L_v^D$ and
\begin{equation}\label{E:26}
(\pp_v,\pp'_v)=\sum_{\delta_i}\frac{|2|_v}{2}\int \L^{D\delta_i}_v(\ppi(\v{a})\pp_v)
\overline{\L^{D\delta_i}_v(\ppi(\v{a})\pp'_v)}\frac{da}{|a|_v}
\end{equation}
is a $G'_v$-invariant Hermitian form. (We used the  factor $\frac{|2|_v}{2}$ to
be consistent with \cite{bm} and \cite{bm2}. There we defined the Hermitian form
on $\widetilde{GL}_2$ first and restricted the form to $\widetilde{SL}_2$.
See section (9.7) of \cite{bm}).
For some explicit constructions of this form, see \S~\ref{S:local}.

We can now define the constant $d_{\ppi}(S,\psi^D)$.  From
the uniqueness of Hermitian forms and Whittaker models, there exist constants
 $\tilde{c}_1(\ppi,S,\psi^D,\{\L_v^D\})$ and 
$\tilde{c}_2(\ppi,S,\psi^D,\{\L_v^D\})$ such that 
whenever $\pp=\otimes_{v\in S}\pp_v\otimes_{v\not\in S}\pp_{0,v}$
(under an identification between
$V_{\ppi}$ and the restricted tensor product $\otimes V_{\ppi,v}$):
$$\W_{\pp}^D(e)=\tilde{c}_1(\ppi,S,\psi^D,\{\L_v^D\})\prod_{v\in S}\L_v^D(\pp_v),$$
$$\|\pp\|=\tilde{c}_2(\ppi,S,\psi^D,\{\L_v^D\})\prod_{v\in S}\|\pp_v\|.$$
As in the case of Lemma~\ref{L:defd1}, we have:

\begin{lemma}\label{L:defd2}
The constant
$$d_{\ppi}(S,\psi^D)=|\tilde{c}_1(\ppi,S,\psi^D,\{\L^D_v\}/
\tilde{c}_2(\ppi,S,\psi^D,\{\L^D_v\})|$$
is independent of the choice of the linear forms $\L^D_v$.
\end{lemma}

When $\ppi$ does not have a nontrivial $\psi^D$-Whittaker model, we
will set $d_{\ppi}(S,\psi^D)=0$. 
The constant
$d_{\ppi}(S,\psi^D)$ is well defined with our fixed choice of $\psi^D$ and 
the measure on $G'$, (and again it is independent of the choice of additive
measure). Explicitly for any vector $\pp=
\otimes_{v\in S}\pp_v\otimes_{v\not\in S}
\pp_{0,v}$ with $\L^D_v({\pp_v})\not=0$ for $v\in S$,
\begin{equation}\label{E:d2explicit}
d_{\ppi}(S,\psi^D)=\frac{|\W_{\pp}^D(e)|}{\|\pp\|}\prod_{v\in S}
\frac{\|\pp_v\|}{|\L^D_v(\pp_v)|}.
\end{equation}

\section{Results on the theta correspondence}\label{S:Waldspurger}

 As stated in the introduction, the generalization of the Shimura correspondence between
 the half integral weight modular forms and integral weight modular forms is the theta 
 correspondence between automorphic representations of $PGL_2({\bf A})$ and $\widetilde{SL}_2
 ({\bf A})$. The theta correspondence is dependent on the choice of the additive
 character $\psi$. 
 With a fixed $\psi$, we denote by $\Theta(\pi,\psi)$
 the representation of $\widetilde{SL}_2({\bf A})$ associated to $\pi$ of 
 $PGL_2({\bf A})$, 
  and $\Theta(\pi_v,\psi_v)$ the representation of  $\widetilde{SL}_2(F_v)$
 associated to $\pi_v$ of $PGL_2(F_v)$ under the 
 theta correspondence. Conversely, the theta correspondence
 associates to $\ppi$ of $\widetilde{SL}_2({\bf A})$ and  $\ppi_v$ 
 of  $\widetilde{SL}_2(F_v)$ representations
 $\Theta(\ppi,\psi)$ of $PGL_2({\bf A})$ and $\pi_v=\Theta(\ppi_v,\psi_v)$ of
 $PGL_2(F_v)$. 

  In this paper, of particular importance are
 the representation $\tilde{\pi}^D=\Theta(\pi\otimes\chi_D,\psi^D)$
 and its local counterpart $\ppi_v^D=\Theta(\pi_v\otimes \chi_D,\psi_v^D)$.
  In this section, we recall Waldspurge's beautiful results on these representations.
  The  works in  \cite{wa},
 \cite{wa2} tell
  us that the set  $\{\ppi^D\}$ (or $\{\ppi_v^D\}$) is finite,
 moreover the dependence of $\ppi^D$ (or $\ppi_v^D$) on $D$ is also given.

 We first recall Waldspurger's local theory. Fix a place $v$ of $F$. Let
 $P_{0,v}$ be the set of  special or supercuspidal representations (or 
 discrete series represenations when $F_v={\bf R}$) of $PGL_2(F_v)$. For
 $D\in F_v^*$, define $(\frac{D}{\pi_v})\in \pm 1$ by:
 $$(\frac{D}{\pi_v})=\chi_D(-1)\epsilon(\pi_v,1/2)/\epsilon(\pi_v\otimes\chi_D, 1/2).$$
 We then get a partition of $F_v^*=F_v^+\cup F_v^-$ where
 $$F_v^{\pm}(\pi_v)=\{D\in F^*_v|\, (\frac{D}{\pi_v})=\pm 1\}.$$
 \begin{theorem}\label{T:Wlocal} \cite{wa2}

 When $\pi_v\not\in P_{0,v}$, $F_v^+=F_v^*$ and $\ppi^D_v=\Theta(\pi_v,\psi_v)$.

 When $\pi_v\in P_{0,v}$, there are two representations $\ppi^+_v$ and
 $\ppi^-_v$ of $G'_v$, such that $\ppi^D=\ppi^+_v=\Theta(\pi_v,\psi_v)$ when
  $D\in F_v^+$, and
 $\ppi_v^D=\ppi^-_v$ when $D\in F_v^-$. 

 Moreover $\ppi^D=\Theta(\pi_v,\psi_v)$ if and only if
  $\Theta(\pi_v,\psi_v)$ has a nontrivial $\psi_v^D$-Whittaker
 model.
 \end{theorem}

 We now state the global counterpart of this Theorem. 
 Let $\tilde{A}_{00}$ be the space of cuspidal automorphic
 forms on $G'(\aa)$ that are orthogonal to the theta series generated by
 quadratic forms of one variable. Let $A_{0,i}$ be the subspace of
 cuspidal automorphic forms on $PGL_2(\aa)$, such that for any $\pi$ subrepresentation
 of $A_{0,i}$, there is $D\in F^*$, with $L(\pi\otimes \chi_D,1/2)\not=0$.

 For $\ppi_1,\ppi_2$ irreducible subrepresentations of $\tilde{A}_{00}$, we
 will say $\ppi_1\sim \ppi_2$ if they are near equivalent, that is at almost
 all places $v$, $\ppi_{1,v}\cong \ppi_{2,v}$. Denote by $\bar{A}_{00}$
 the quotient of $\tilde{A}_{00}$ by this relation. 

 Let $\Sigma=\Sigma(\pi)$ be the set of places $v$ where
 $\pi_v\in P_{0,v}$.
 Given $D\in F^*$, let $\epsilon(D,\pi)=(\frac{D_v}{\pi_v})_{v\in \Sigma}$.
 Then $\epsilon(D,\pi)\in \{\pm 1\}^{|\Sigma|}$. We have
 \begin{equation}\label{E:p3}
 \epsilon(\pi\otimes\chi_D,1/2)=\epsilon(\pi,1/2)\prod_{v\in\Sigma}(\frac{D_v}{\pi_v}).
 \end{equation}
 We will use $\epsilon=
 (\epsilon_v)_{v\in \Sigma}$ to denote an element in $\{\pm 1\}^{|\Sigma|}$,
 with $\epsilon_v\in\{\pm 1\}$. Given such an $\epsilon$, we will let $F^{\epsilon}(\pi)$ to
 be the set of $D\in F^*$ with $\epsilon(D,\pi)=\epsilon$. Then we get a partition
 $F^*=\cup_{\epsilon\in \{\pm 1\}^{|\Sigma|}}F^{\epsilon}(\pi)$.

 \begin{theorem}\label{T:Wglobal} \cite{wa2}

 1. (Relation with local correspondence) When $\Theta(\ppi,\psi)\not=0$, $\Theta(\ppi,\psi)\cong\otimes_v\Theta(\ppi_v,\psi_v)$.
 When $\Theta(\pi,\psi)\not=0$, $\Theta(\pi,\psi)\cong\otimes_v\Theta(\pi_v,\psi_v)$.

 2. (Nonvanishing of the correspondence)
 $\Theta(\pi,\psi)\not=0$ if and only if $L(\pi,1/2)\not=0$. $\Theta
 (\ppi,\psi)\not=0$ if and only if $\ppi$ has a nontrivial $\psi$-Whittaker
 model.

 3. (Correspondence as a bijection)
 For $\ppi$ an irreducible subrepresentation of $\tilde{A}_{00}$, there is a unique
 $\pi$ associated to $\ppi$, such that whenever $\Theta(\ppi,\psi^D)\not=0$,
 $\Theta(\ppi,\psi^D)\otimes \chi_D=\pi$. Denote this association by
 $\pi=S_{\psi}(\ppi)$.
 This association defines a bijection between $\bar{A}_{00}$ and
 $A_{0,i}$.

 4. (Description of near equivalent class).
 If $\pi=S_{\psi}(\ppi)$,
 the near equivalence class of $\ppi$ consists of all the  nonzero $\ppi^D$'s. 

 5. (Dependence of $\ppi^D$ on $D$). Let $\epsilon\in \{\pm 1\}^{|\Sigma|}$.
 If  $\prod_{v\in \Sigma}\epsilon_v\not=\epsilon(\pi,1/2)$, then
 $\ppi^D=0$ for all $D\in F^{\epsilon}(\pi)$. If $\prod_{v\in \Sigma}\epsilon_v=\epsilon(\pi,1/2)$,
 then there is a unique $\ppi^{\epsilon}$ such that
 for $D\in F^{\epsilon}(\pi)$,  $\ppi^D=\ppi^{\epsilon}$ when
  $L(\pi\otimes\chi_D,1/2)\not=0$ and $\ppi^D=0$ otherwise.
 \end{theorem}

 For convenience, if $\prod_{v\in \Sigma}\epsilon_v\not=\epsilon(\pi,1/2)$, we set $\ppi^{\epsilon}=0$.

\section{Statement of the main results}\label{S:mainresult}

\subsection{The formula for $L(\pi,1/2)$}

The definition of the $L-$function $L(\pi,s)$ and the local $L-$functions $L(\pi_v,s)$
can be found in \cite{jl}.
Fix a finite set of places $S$, we use $L^S(\pi,s)$ to denote
 the partial $L-$function $\prod_{v\not\in S}L(\pi_v,s)$.
\begin{theorem}\label{T:21}
For an irreducible cuspidal
automorphic representations $\pi$ of $GL_2(\aa)$ with trivial
central character and  $L(\pi,1/2)\not=0$, for 
$D\in F^{*}$, let $\ppi_D=\Theta(\pi,\psi^D)$. Let
 $S$ be a finite set of places as in the introduction and moreover containing all places of $v$
where $\psi$ or $\psi^D$ is not unramifed, and all places  where $\pi_v$ or
$\ppi_{D,v}$  is not unramified. Then
\begin{equation}\label{E:28}
|d_{\pi}(S,\psi)|^2L^S(\pi,1/2)=|d_{\ppi_D}(S,\psi^D)|^2
\end{equation}
\end{theorem}

We state a more explicit formula using (\ref{E:d1explicit}) and
(\ref{E:d2explicit}). Take any vectors
$\varphi=\otimes_{v\in S}\varphi_v\otimes_{v\not\in S}
\varphi_{0,v}$ and $\pp=\otimes_{v\in S}\pp_v\otimes_{v\not\in S}\pp_{0,v}$ 
in $V_{\pi}$ and $V_{\ppi_D}$ such that 
$L_v(\varphi_v)\not=0$ and $\L_v(\pp_v)\not=0$.
Define:
\begin{eqnarray}
e(\varphi_v,\psi_v)&=&\frac{\|\varphi_v\|^2}{|L_v(\varphi_v)|^2}\\
e(\pp_v,\psi_v^D)&=&\frac{\|\pp_v\|^2}{|\L^D_v(\pp_v)|^2}
\end{eqnarray}
Then as in the proof of the Lemma~\ref{L:defd1}, these constants are
independent of our choice of $L_v$ and $\L^D_v$ and are well defined. 
From (\ref{E:d1explicit}) and
(\ref{E:d2explicit}), we see the identity (\ref{E:28}) can be stated as follows:
\begin{equation}\label{E:main1}
\frac{|\tilde{W}_{\pp}^D(e)|^2}{\|\pp\|^2}=\frac{|W_{\varphi}(e)|^2L(\pi,1/2)}{
\|\varphi\|^2}\prod_{v\in S}\frac{e(\varphi_v,\psi_v)}{e(\pp_v,\psi_v^D)L(\pi_v,1/2)}
\end{equation}

\subsection{The adelic version of Waldspurger's Theorem}

The statement in Theorem~\ref{T:21} is a direct result of our interpretation of Waldspurger's
formula. To see the relation with the previous versions of Waldspurger's formula,
we apply the theorem to the case with $\pi$ replaced by $\pi\otimes\chi_D$. Then 
$\ppi_D=\ppi^D=\Theta(\pi\otimes\chi_D,\psi^D)$ and we can apply the results in \S~\ref{S:Waldspurger}.
 The results are stated in the
following two Theorems.

In the rest of the section, we will assume $D$ satisfies that
 for all {\em odd} nonarchimedean places $v$, $|D|_v=1$ or $|D|_v=q_v^{-1}$, and for all
even nonarchimedean places $v$, $1\leq |D|_v\leq q_v^{-2}$.
With a bit abuse of terminology, we call this $D$
 a square free integer
in $F^*$. 

In \cite{wa1}, Waldspurger described the Fourier coefficient of a half integral weight form
using the data from the corresponding integral weight form. We generalize his result here.
For 
$\ppi$ an irreducible subrepresenation  of $\tilde{A}_{00}$
and $D\in F^*$, we describe $d_{\ppi}(S,\psi^D)$ 
in terms of  the data of $\pi=S_{\psi}(\ppi)$.

Let $\Sigma=\Sigma(\pi)$ be
as in Theorem~\ref{T:Wglobal}. From Theorem~\ref{T:Wglobal},
$\ppi=\ppi^{\epsilon_0}$ for some $\epsilon_0\in \{\pm 1\}^{|\Sigma|}$.

\begin{theorem}\label{T:Wadele}
Let $S$ be as in the introduction and contain the places
$v$ where $\ppi_v$ or $\psi$ is not unramified.
 Let $D$ be a square free integer in $F^*$. Let $\pi=S_{\psi}(\ppi)$ and
$\epsilon_0$ be as above.
If $D\in F^{\epsilon_0}(\pi)$ then 
\begin{equation}\label{E:Wadele}
|d_{\ppi}(S,\psi^D)|^2=|d_{\pi}(S,\psi)|^2L^S(\pi\otimes\chi_D,1/2)\prod_{v\in S}|D|^{-1}_v.
\end{equation}
If $D\not\in F^{\epsilon_0}(\pi)$, $d_{\ppi}(S,\psi^D)=0$.
\end{theorem}

\subsection{Adelic version of the Kohnen-Zagier formula}

The Kohnen-Zagier formula quoted in the introduction, as opposed to the Waldspurger's
theorem in \cite{wa1}, describes the twisted $L-$value of $f(z)$ in terms of data of a
half integral weight form. We now state its generalization.

Let $\pi\in A_{0,i}$. Let $\Sigma=\Sigma(\pi)$ be as in Theorem~\ref{T:Wglobal}.
For $\epsilon\in \{\pm 1\}^{|\Sigma|}$, we define $\ppi^{\epsilon}$
as in Theorem~\ref{T:Wglobal}.

\begin{theorem}\label{T:KZadele}
Let $S$ be as in the introduction
and contain the places where $\psi$  or $\pi_v$ is
not unramified. 
Then for $\epsilon\in  \{\pm 1\}^{|\Sigma|}$, for
$D\in F^{\epsilon}(\pi)$ a square free integer, 
\begin{equation}\label{E:KZadele}
|d_{\ppi^{\epsilon}}(S,\psi^D)|^2=|d_{\pi}(S,\psi)|^2L^S(\pi\otimes\chi_D,1/2)\prod_{v\in S}|D|^{-1}_v.
\end{equation}
\end{theorem}

\begin{remark}
1. The equation (\ref{E:KZadele}) is in fact a finite set of equations, corresponding
to the finite partition of $F^*$ by $F^{\epsilon}(\pi)$.
The conditions on $D$ in (\ref{E:11}) are precisely the condition
  $D\in D^{\epsilon_0}(\pi)$ for a given $\epsilon_0$. Thus (\ref{E:11}) is
only one in a set of equations. See \S~\ref{S:KZ} for the whole set of equations.

2. From Theorem~\ref{T:Wadele}, $|d_{\ppi^{\epsilon}}(S,\psi^D)|=0$
when $D\not\in F^{\epsilon}(\pi)$. Thus for $\pi\in A_{0,i}$, for all $D$
square free integers,
\begin{equation}\label{E:sumKZ}
\sum_{\epsilon}
|d_{\ppi^{\epsilon}}(S,\psi^D)|^2=|d_{\pi}(S,\psi)|^2L^S(\pi\otimes\chi_D,1/2)\prod_{v\in S}|D|^{-1}_v.
\end{equation}

3.
If $\pi\not\in A_{0,i}$, then clearly $L(\pi\otimes \chi_D,1/2)=0$
for all $D\in F^*$.
\end{remark}

\subsection{Some other consequences}

Since $L(\pi_v,1/2)>0$ when $\pi_v$ is unitary, from the Theorem~\ref{T:21} we have
\begin{corollary}
For all irreducible cuspidal
automorphic representations $\pi$ of $PGL_2(\aa)$, $L(\pi, 1/2)\geq 0$.
\end{corollary}
 This result is first
shown in [Gu].

The next corollary can be considered as the adelic version of corollary~2 in \cite{wa1}.
It follows immediately from Theorem~\ref{T:KZadele}.
\begin{corollary}
Let $\pi\in A_{0,i}$ and $S$, $\Sigma=\Sigma(\pi)$ be as in Theorem~\ref{T:KZadele}.
Fix $\epsilon\in \{\pm 1\}^{|\Sigma|}$.
For two square free integers $D_1,D_2\in F^{\epsilon}(\pi)$,
\begin{equation}\label{E:2D}
|d_{\ppi^{\epsilon}}(S,\psi^{D_1})|^2L(\pi\otimes\chi_{D_2},1/2)
=|d_{\ppi^{\epsilon}}(S,\psi^{D_2})|^2L(\pi\otimes\chi_{D_1},1/2)\prod_{v\in S}|\frac{D_2}{D_1}|_v.
\end{equation}
\end{corollary}

We also state a result concerning the relation
of the $D-$th and $D\Delta^2$-th Fourier coefficients of a half-integral
weight form. Note the similarity with Lemma~\ref{L:dpi}.
\begin{corollary}\label{C:delta2}
Let $S$ be as in Theorem~\ref{T:KZadele}.
If $D'=D\Delta^2$ for $\Delta\in F^{\times}$. Let $S_{D,D'}$ be a
finite set of places with $|D|_v=|D'|_v=1$ 
for all $v\not\in S_{D,D'}$. Then
\begin{equation}\label{E:delta2}
|d_{\ppi}(S_{D,D'}\cup S,\psi^D)|=|d_{\ppi}(S_{D,D'}\cup S,\psi^{D'})|.
\end{equation}
\end{corollary}

\subsection{Ramanujuan conjecture and Lindel\H{o}f Hypothesis}

The Ramanujuan conjecture for the half integral weight cusp form is as follows:

\begin{quote}
Let $g(z)=\sum_{n=1}^{\infty}c(n)e^{2\pi inz}$ be a cusp form of 
weight $k+1/2$, $k\in {\bf Z}$,
such that $g(z)$ is orthogonal to the space generated by
the theta series associated to quadratic forms of one variable, (i.e.
$g(z)$ is a vector in $\tilde{A}_{00}$). Then when $n$ is square free, as
$n\mapsto \infty$,
\begin{equation}\label{E:rc1}
|c(n)|<<_{g,\alpha}n^{k/2-1/4+\alpha}
\end{equation}
for all $\alpha>0$. The implied constant depends on $g(z)$ and $\alpha$ only.
\end{quote}

To state the generalization of this conjecture to the case of cusp forms
over totally real fields, we find it most natural to use the notion $d_{\ppi}(S,\psi^D)$
again.

Let $S$ be a finite set of places as in the introduction and contain all places where $\ppi_v$
is not unramified. 
For $S_0\subset S$,  for $\pp$ a vector in the space of $\ppi$,
we define
\begin{equation}\label{E:d2variant}
d_{\ppi}(\pp, S_0,\psi^D)=\frac{|\W_{\pp}^D(e)|}{\|\pp\|}\prod_{v\in S_0}
\frac{\|\pp_v\|}{|\L^D_v(\pp_v)|}=d_{\ppi}(S,\psi^D)\prod_{v\in S-S_0}
\frac{|\L^D_v(\pp_v)|}{\|\pp_v\|}.
\end{equation}
As before this constant is well defined and independent of the choice of $\{\L^D_v\}$.

Let $S_{\infty}$ be the collection of archimedean places of $F$. For $D\in F^*$, define
$|D|_{\infty}=\prod_{v\in S_{\infty}}|D|_v$.
\begin{conjecture}\label{rc}
(Ramanujuan conjecture).
Let $\ppi$ be an irreducible subrepresentation of $\tilde{A}_{00}$. Let
$\pp\in V_{\ppi}$. 
For $D$ a square free integer in $F^*$, as $|D|\mapsto \infty$, for all $\alpha>0$
\begin{equation}\label{E:rc2}
|d_{\ppi}(\pp,S_{\infty},\psi^D)|<<_{\ppi,\pp,\alpha}|D|_{S_{\infty}}^{\alpha-1/2}
\end{equation}
where the implied constant depends only on $\ppi,\pp$ and $\alpha$.
\end{conjecture}

We will show in \S~\ref{S:dictionary} that the above conjecture implies inequality  (\ref{E:rc1}).

The Lindel\H{o}f hypothesis is a conjecture on the bound of central value
of $L-$functions. We state only a special case.
\begin{conjecture}\label{lh}
(Lindel\H{o}f hypothesis) Let $\pi$ be an irreducible cuspidal automorphic representation
of $G(\aa)$ with trivial central character, then for 
$D$ square free integer, as $|D|\mapsto \infty$, for all $\beta>0$
\begin{equation}\label{E:lh}
|L^{S_{\infty}}(\pi\otimes \chi_D,1/2)|<<_{\pi,\beta}|D|_{S_{\infty}}^{\beta}
\end{equation}
where the implied constant depends only on $\pi$ and $\beta$.
\end{conjecture}
\begin{theorem}\label{T:conjecture}
The inequality (\ref{E:rc2}) holds for some $\alpha>0$ (and for all $\ppi$ and $\pp$ as
in Conjecture~\ref{rc}) if and only if the inequality (\ref{E:lh}) holds for 
$\beta=2\alpha>0$ (for all $\pi$ as in Conjecture~\ref{lh}). In particular, the Conjecture~\ref{rc}
is equivalent to the Conjecture~\ref{lh}.
\end{theorem}

\section{A relative trace formula}\label{S:Jacquet}

In \cite{j}, Jacquet proved some of  Waldspurger's results
on theta correspondence using a relative trace
formula. Our result is based on a local analysis of his trace
formula and its variation. We recall some of the results on the trace
formula, in the process fix some notations. The main result
here is Theorem~\ref{T:35}.

\subsection{Definition of the global distributions $I(f,\psi)$ and $J(f',\psi^D)$}

 Let $f(g)\in \cc(
Z({\bf A})\backslash
G(\bf A))$ the space of smooth compactly supported functions. Define a kernel function
\[K_f(x,y)=\sum_{\xi\in PGL_2(F)}f(x^{-1}\xi y)\]
Define a distribution $I(f,\psi)$ to be:
\begin{equation}\label{E:31}
\int_{{\bf A}^{*}/F^{*}}\int_{{\bf A}/F}
K_f(\u{a},n(u))\psi(u)du
d^{*}a\end{equation}
Let $f'(g)\in \cc(G'(\bf A))$, (we use this notation to denote the space of
genuine smooth compactly supported functions). Define a kernel function
\[K_{f'}(x,y)=\sum_{\xi\in SL_2(F)}f'(x^{-1}\cdot\xi\cdot y)\]
Here we note that $SL_2(F)$ embeds into $G'(\bf A)$.
Define a distribution $J(f',\psi^D)$ to be:
\begin{equation}\label{E:32}
\int_{{\bf A}/F}\int_{{\bf A}/F}K_f(\tilde{n}(x),\tilde{n}(y))\psi^D
(-x+y)dxdy\end{equation}

The relative trace formula is an identity between the distributions
$I(f,\psi)$ and $J(f',\psi^D)$. 
We will state the result on the distributions $I(f,\psi)$, $J(f',
\psi^D)$. The computations are available in \cite{j} and will not be
included.

\subsection{Comparison of orbital integrals}

\begin{proposition}\label{P:31} \cite{j}
If $f=\otimes f_v$ and $f'=\otimes f'_v$, then
\[I(f,\psi)=\prod I_{\psi}^+(f_v)+\prod I_{\psi}^-(f_v)+\sum_{a\in F^{*}}\prod
{\mathcal O}^{\psi_v}_{f_v}(n(a)w)\]
\[J(f',\psi^D)=\prod J_{\psi^D}^+(f'_v)+\prod J_{\psi^D}^-(f'_v)+\sum_{a\in F^{*}}
\prod {\mathcal O}^{\psi_v^D}_{f'_v}(\tilde{w}\cdot \v{a})\]
\end{proposition}
In the above equations, $I_{\psi}^{\pm}(f_v)$ and $J_{\psi^D}^{\pm}(f'_v)$
are the so called
{\em singular orbital integrals} of $f_v$ and $f'_v$, whose precise
forms are not important for us,
while $${\mathcal O}^{\psi_v}_{f_v}(g)=\int f_v(\u{a}gn(x))\psi_v(x)dxd^{*}a$$
$${\mathcal O}^{\psi_v^D}_{f'_v}(g)=\int f'_v(\tilde{n}(x)\cdot g\cdot\tilde{n}(y))
\psi_v^D(x+y)dxdy$$
\begin{proposition}\label{P:32} \cite{j}
For each $f$ in $\cc(G(F_v)/Z(F_v))$, there is
$f'\in \cc(G'(F_v))$
such that for $a\in F_v^{*}$
\begin{equation}\label{E:34}
{\mathcal O}^{\psi_v}_{f}(n(\frac{a}{4D})w)=
{\mathcal O}^{\psi_v^D}_{f'}(\tilde{w}\cdot \v{a})\psi_v(-\frac{2D}{a})
|a|_v^{1/2}\gamma(a^{-1},\psi_v^D)^{-1}
\end{equation}
and $I^{\pm}_{\psi}(f)=J^{\pm}_{\psi^D}(f')$.
Conversely, given $f'$, we can find a $f$
satisfying the equations.
\end{proposition}
 We say the two functions $f$ and $f'$ match if the relations
in the proposition are satisfied.

Now let $v$ be a non-Archimedean place with odd residue characteristic,
and where $\psi_v$, $\psi_v^D$ have order 0. Recall
that the Hecke algebra ${\mathcal H}(G(F_v)/Z(F_v))$
is the algebra of compactly supported functions
on $G(F_v)/Z(F_v)$ biinvariant under the maximal compact group $G(\oo_v)$. The
Hecke algebra ${\mathcal H}(G'(F_v))$
is similarly defined, except that the functions
are genuine, and biinvariant under $SL_2(\oo_v)$ embedded in
$G'(F_v)$.
\begin{proposition}\label{P:33} \cite{j}
There is an algebra isomorphism $\eta_v: {\mathcal H}(G(F_v)/Z(F_v))\rightarrow
{\mathcal H}(G'(F_v))$, such that $f$ and $\eta_v(f)$ match.
\end{proposition}

From Propositions~\ref{P:31}, \ref{P:32} and \ref{P:33},
 we get
\begin{theorem}\label{T:34}
Fix any finite set of places $S$ as in the introduction and containing
places where $\psi$ and $D$ are not of order 0. For each place $v\in S$, there is a map
$\rho_v: \cc(G(F_v)/Z(F_v))\rightarrow
\cc(G'(F_v))
$, such that 
$$I(\otimes_{v\in S} f_v\otimes_{v\not\in S}f_v,\psi)=J(\otimes_{v\in S}\rho_v(f_v)
\otimes_{v\not\in S}\eta_v(f_v),\psi^D),$$ whenever 
$f_v\in {\mathcal H}(G(F_v)/Z(F_v))$ for all $v\in S$.
\end{theorem}

\subsection{Relation with Shimura-Waldspurger correspondence}

For $\pi$ an irreducible cuspidal automorphic representation of $G(\aa)$
with trivial central character,
define
\begin{equation}\label{E:35}
I_{\pi}(f,\psi)=\sum_{\varphi_i}Z(\pi(f)\varphi_i)
\overline{{W}_{\varphi_i}(e)}\end{equation}
with $\varphi_i$ an orthonormal basis of $V_{\pi}$;
here for $\varphi\in V_{\pi}$
$$\pi(f)\varphi=\int_{Z(\aa)\backslash G(\aa)}f(g)\pi(g)\varphi dg,$$
$$Z(\varphi)=\int_{\aa^*/F^*}\varphi(\u{a})d^{*}a.$$
For
$\ppi$ an irreducible cuspidal representation of $G'(\aa)$, define
\begin{equation}\label{E:36}
J_{\ppi}(f',\psi^D)=\sum_{\varphi'_i}\W^D_{\ppi(f')\pp_i}(e)
\overline{{\W}^D_{\pp_i}(e)}\end{equation}
with $\varphi'_i$ an orthonormal basis of $V_{\ppi}$; here for $\pp\in V_{\ppi}$
$$\ppi(f')\pp=\int_{ G'(\aa)}f'(g)\ppi(g)\pp dg.$$

The distributions $I_{\pi}(f,\psi)$ and $J_{\ppi}(f',\psi^D)$ are
the contributions from $\pi$ and $\ppi$ to $I(f,\psi)$ and $J(f',\psi^D)$
respectively.

Recall that if $\pi_v$ is unramified, there exists a vector
$\varphi_{0,v}$ that is fixed under the action of $G(\oo_v)$.
For $f_v$ a Hecke
function on $G(F_v)/Z(F_v)$, there is a constant $\hat{f}_v(\pi_v)$ with
\begin{equation}\label{E:37}
\pi_v(f_v)\varphi_{0,v}=\hat{f}_v(\pi_v)\varphi_{0,v}.
\end{equation}
Similarly, if $\ppi_v$ be unramified, let $\pp_{0,v}$
be a vector that is fixed under $SL_2(\oo_v)$, then
for $f'_v$ in the Hecke algebra of $G'(F_v)$, there is
a constant $\hat{f'}_v(\ppi_v)$ with \begin{equation}\label{E:38}
\ppi_v(f'_v)\pp_{0,v}
=\hat{f'}_v(\ppi_v)\pp_{0,v}.\end{equation}

It is standard to derive from the Theorem~\ref{T:34} the following:

\begin{theorem}\label{T:35} \cite{j}
For any cuspidal representation $\pi$ of $G$ with trivial
central character such that
$I_{\pi}(f,\psi)$
is nontrivial, there is a unique cuspidal representation $\ppi$ of
$G'$, such that if $f$ and $f'$ match
\begin{equation}
I_{\pi}(f,\psi)=J_{\ppi}(f',\psi^D)
\end{equation}
Moreover, if $S$ satisfies the condition in Theorem~\ref{T:34} and
contains all places where $\pi_v$ or $\ppi_v$ is not unramified,
for $v\not\in S$, if $f_v$ is a Hecke function and $f'_v=\eta_v(f_v)$,
then 
\begin{equation}\label{E:satake}
\hat{f}_v(\pi_v)=\hat{f'}_v(\ppi_v).
\end{equation}
\end{theorem}
\begin{remark}
From the definition of $Z(\varphi)$ and the integral representation of $L-$function
$L(\pi,s)$, it is clear that $I_{\pi}(f,\psi)$ is nontrivial if and
only if  $L(\pi,1/2)\not=0$.
\end{remark}

\begin{proposition}\label{P:theta}
In the above theorem, $ \ppi=\Theta(\pi,\psi^D)$.
\end{proposition}
\begin{proof}
From the description of the map $\eta_v$ of Hecke algebras in \cite{j}
and the equation (\ref{E:satake}), we get $\ppi$ is in the same
near equivalence class as $\Theta(\pi,\psi^D)$. As $J_{\ppi}(f',\psi^D)\not=0$,
$\ppi$ has $\psi^D-$Whittaker model from (\ref{E:36}). Let $\Sigma=\Sigma(\pi)$ be
 as in Theorem~\ref{T:Wglobal}. Then for $v\in \Sigma$,
$\ppi_v=\Theta(\pi_v,\psi_v^D)$ by Theorem~\ref{T:Wlocal}. Thus $\ppi$ must
be the same representation as $\Theta(\pi,\psi^D)$. 
\end{proof}

\section{The local distributions}\label{S:BM}

Let $\pi$ and $\ppi$ be the cuspidal representations that correspond
to each other by Theorem~\ref{T:35}. Then $\ppi=\Theta(\pi,\psi^D)$.
Let $S$ be as in Theorem~\ref{T:35}.
Assume $f=\otimes f_v, f'=\otimes f'_v$, where $f$ and $f'$
match, and $f_v$, $f'_v$ are matching Hecke functions
when $v\not\in S$. We write
$I_{\pi}(f,\psi)$ and $J_{\ppi}(f',\psi^D)$ as products of
local distributions over the places in $S$.
We then state the identity between the local distributions,
which coupled with Theorem~\ref{T:35} gives Theorem~\ref{T:21}.

\subsection{The distribution $I_{\pi,v}(f_v,\psi_v)$}

We fix a choice of local Whittaker functionals $L_v$ on $\pi_v$,
and define the Hermitian form on $V_{\pi_v}$ using this choice of
$L_v$. 
For $v\in S$, we fix for $\pi_v$ as above
an orthonormal basis of $V_{\pi,v}$, denote it by $\{\varphi_{i,v}\}$.
For $v\not\in S$, let $\varphi_{0,v}$ be the vector given in \S~\ref{S:constants}. 
For  $\pi=\otimes \pi_v$, from (\ref{E:defc2}), the set
\begin{equation}\label{E:basis}
\{\varphi_I\}=
\{c_2(\pi,S,\psi,\{L_v\})^{-1}\otimes_{v\in S} \varphi_{i,v}\otimes_{v\not\in S}
\varphi_{0,v}\}
\end{equation}
 can be extended to an orthogonal basis of $V_{\pi}$.
Let $V(\pi,S)$ be the space of vectors generated by the set
of $\{\varphi_I\}$. With our choice of $f$,
if is clear that if $\varphi\in V_{\pi}$ is perpendicular to
the space $V(\pi,S)$, then $\pi(f)\varphi=0$.
Thus the expression (\ref{E:35}) for $I_{\pi}(f,\psi)$ takes the form:
\begin{equation}\label{E:40}
\sum_{\varphi_I} Z(\pi(f)\varphi_I)
\overline{W_{\varphi_I}(e)}.
\end{equation}

Using the Hecke theory for $GL_2$, we show:
\begin{lemma}\label{L:Hecke}
When $
\varphi=\otimes_{v\in S}\varphi_{v}\otimes_{v\not\in S}
\varphi_{0,v}$, 
\begin{equation}
Z(\varphi)=c_1(\pi,S,\psi,\{L_v\})
L(\pi,1/2)\prod_{v\in S}\lambda_v(\varphi_{v})
\end{equation}
where
$$\lambda_v(\varphi_{v})=\frac{\int_{F^*_v} L_v(\pi_v(\u{a}){\varphi_{v}})
|a|_v^{s-1/2}d^{*}a}{L(\pi,s)}|_{s=1/2}.$$
\end{lemma}
\begin{proof}
Since $\varphi(\u{a})=\sum_{\delta\in F^*}W_{\varphi}(\u{\delta}\u{a})$, we get:
$$Z(\varphi)=
\int_{\aa^{*}}W_{\varphi}(\u{a})|a|_v^{s-1/2}
d^{*}a|_{s=1/2},$$
which by (\ref{E:defc1}) equals
$$L(\pi,1/2)c_1(\pi,S,\psi,\{L_v\})\prod_{v}\frac{\int_{F^*_v} L_v(\pi_v(\u{a}){\varphi_{v}})
|a|_v^{s-1/2}d^{*}a}{L(\pi,s)}|_{s=1/2}.$$
When $v\not\in S$, it is well known that the above local factor equals 1, (\cite{g}).
Thus the Lemma.
\end{proof}

\begin{proposition}\label{P:41}
Let $S$ be as in Theorem~\ref{T:35}. When $f=\otimes f_v$ where
$f_v$ is a Hecke function if $v\not\in S$:
\begin{equation}\label{E:42}
I_{\pi}(f,\psi)=L(\pi,1/2)|d_{\pi}(S,\psi)|^2\prod_{v\in S}I_{\pi,v}(f_v,\psi_v)
\prod_{v\not\in S}\hat{f}_v(\pi_v)
\end{equation}
where
\begin{equation}\label{E:43}
I_{\pi,v}(f_v,\psi_v)=\sum_{\varphi_{i,v}}\lambda_v(\pi_v(f_v)\varphi_{i,v})
\overline{L_v(\varphi_v)}.\end{equation}
where the sum is taken over the orthonormal basis of $V_{\pi,v}$.
\end{proposition}
\begin{proof}
Let $\varphi=c_2(\pi,S,\psi,\{L_v\})^{-1}\otimes_{v\in S} \varphi_{v}
\otimes_{v\not\in S}
\varphi_{0,v}$
be an element in the orthonormal set (\ref{E:basis}).
From (\ref{E:defc1}), we get 
$$W_{\varphi}(e)=\frac{c_1(\pi,S,\psi,\{L_v\})}{c_2(\pi,S,\psi,\{L_v\})}\prod_{v\in S}L_v(\varphi_v).$$
From (\ref{E:37}),
$$\pi(f)\varphi=c_2(\pi,S,\psi,\{L_v\})^{-1}\prod_{v\not\in S}\hat{f}_v(\pi_v)
\otimes_{v\in S}
\pi_v(f_v)\varphi_v\otimes_{v\not\in S}\varphi_{0,v}.$$
From the above Lemma:
$$Z(\pi(f)\varphi)
\overline{W_{\varphi}(e)}=|\frac{c_1(\pi,S,\psi,\{L_v\})}{c_2(\pi,S,\psi,\{L_v\})}|^2
L(\pi,1/2)\prod_{v\not\in S}\hat{f}_v(\pi_v)
\prod_{v\in S}\lambda_v(\pi_v(f_v)\varphi_{v})\overline{L_v(\varphi_v)}.
$$
The proposition follows from (\ref{E:40}) and the definition of $d_{\pi}(S,\psi)$.
\end{proof}

\begin{remark}
 The expression $I_{\pi,v}(f_v,\psi_v)$ is well defined and independent of the
linear form $L_v$ we choose, as a change in $L_v$ will result in
a change in the Hermitian form, thus the orthonormal basis of $V_{\pi_v}$,
leaving $I_{\pi,v}(f_v,\psi_v)$ unchanged.
\end{remark}

\subsection{The distribution $J_{\ppi,v}(f'_v,\psi_v^D)$}

We can apply the above argument also to $J_{\ppi}(f',\psi_v^D)$. Similarly
we have:

\begin{proposition}\label{P:42}
Let $S$ be as in Theorem~\ref{T:35}. When $f'=\otimes f'_v$ where
$f'_v$ is a Hecke function if $v\not\in S$:
\begin{equation}\label{E:44}
J_{\ppi}(f',\psi^D)=|d_{\ppi}(S,\psi^D)|^2\prod_{v\in S}J_{\ppi,v}(f'_v,\psi_v^D)
\prod_{v\not\in S}\hat{f'}_v(\ppi_v)\end{equation}
where
\begin{equation}\label{E:45}
J_{\ppi,v}(f'_v,\psi_v^D)=\sum_{\pp_{j,v}}\L^D_v(\ppi_v(f'_v)\pp_{j,v})
\overline{\L^D_v(\pp_{j,v})}\end{equation}
where the sum is taken over the orthonormal basis $\{\pp_{j,v}\}$
of $V_{\ppi,v}$.
\end{proposition}

\begin{remark}
 Again one can show that the
 expression $J_{\ppi,v}(f'_v,\psi_v^D)$ is well defined and independent of the
linear form $\L^D_v$ we choose.
\end{remark}

\subsection{Statement of the local identity}

\begin{theorem}\label{T:43}
Fix a place $v$, when $f_v,f'_v$ match, when $\pi_v$ is
a local component of an irreducible cuspidal automorphic representation
$\pi$ of $PGL_2(\aa)$ with $L(\pi,1/2)\not=0$,  let
$\ppi_v=\Theta(\pi_v,\psi_v^D)$,
then \begin{equation}\label{E:46}
J_{\ppi,v}(f'_v,\psi_v^D)=|2D|_v\epsilon(\pi_v,1/2)
L(\pi_v,1/2)I_{\pi,v}(f_v,\psi_v)\end{equation}
\end{theorem}

The proof of this Theorem is quite technical. It is done in \cite{bm},\cite{bm2}.
In \cite{bm} we established the identity when $v$ is nonarchimedean. In
\cite{bm2}, the identity is proved when $v={\bf R}$. In this case, the
identity follows from the identities between classical Bessel functions. We
established the Theorem in a bit more generality, as we do not
assume $\pi_v$ is a local component of an automorphic representation.
 The proof of the Theorem stated
as above is easier, as from Theorem~\ref{T:35}, we have $
J_{\ppi,v}(f'_v,\psi_v^D)=c(\pi_v,\ppi_v)I_{\pi,v}(f_v,\psi_v)$ for
some constant $c(\pi_v,\ppi_v)$
independent of $f_v$ and $f'_v$. One then only needs to
determine this constant.

\section{Proof of the main results}\label{S:proof}

We now prove the results stated in \S~\ref{S:mainresult}.

{\sc Proof of Theorem~\ref{T:21}}: Let $\pi$ and $\ppi=\ppi_D$
be as in the Theorem. From Theorem~\ref{T:35} and Proposition~\ref{P:theta},
we see $I_{\pi}(f,\psi)=J_{\ppi}(f',\psi^D)$ when $f$ and $f'$ 
match. Assume $f=\otimes f_v$ with $f_v$ a Hecke function if $v\not\in S$,
and $f'=\otimes f'_v$ with $f'_v=\eta_v(f_v)$ when $v\not\in S$ and
$f'_v$ matches $f_v$ elsewhere. From Theorem~\ref{T:35} and Propositions
\ref{P:41} and \ref{P:42}, we get:
\begin{equation}\label{E:proof1}
L(\pi,1/2)|d_{\pi}(S,\psi)|^2\prod_{v\in S}I_{\pi,v}(f_v,\psi_v)
=|d_{\ppi}(S,\psi^D)|^2\prod_{v\in S}J_{\ppi,v}(f'_v,\psi_v^D).
\end{equation}
From Theorem~\ref{T:43}, we get:
\begin{equation}\label{E:proof2}
\prod_{v\in S}J_{\ppi,v}(f'_v,\psi_v^D)
=\prod_{v\in S}|2D|_v\epsilon(\pi_v,1/2)
L(\pi_v,1/2)I_{\pi,v}(f_v,\psi_v).
\end{equation}
Combine the above two equations, we get:
$$L(\pi,1/2)|d_{\pi}(S,\psi)|^2=|d_{\ppi}(S,\psi^D)|^2
\prod_{v\in S}|2D|_v\epsilon(\pi_v,1/2)
L(\pi_v,1/2).$$
As $|2D|_v=1$ for $v\not\in S$, we get $\prod_{v\in S}|2D|_v=1$.
As $\epsilon(\pi,1/2)=1$ and for $v\not\in S$, $\epsilon(\pi_v,1/2)=1$,
we get $\prod_{v\in S}\epsilon(\pi_v,1/2)=1$. Thus we get the identity
(\ref{E:28}). $\qed$

{\sc Proof of Theorem~\ref{T:Wadele}:} 
First assume $D\in F^{\epsilon_0}(\pi)$. Let $S_D$ be a finite set of
places, such that $|D|_v=1$ when $v\not\in S_D$. Let $S_1=S\cup S_D$,
$S_2=S_1-S$. Let $\pi=S_{\psi}(\ppi)$ be as in the Theorem and
$\ppi^D=\Theta(\pi\otimes \chi_D,\psi^D)$. 
Then by Theorem~\ref{T:Wglobal},
we have $\ppi^D=\ppi$ or 0. 

When $\ppi^D=0$, $L(\pi\otimes\chi_D,1/2)=0$
from Theorem~\ref{T:Wglobal}. Meanwhile from Propsition~30 of \cite{wa},
$\ppi$ does not have a $\psi^D-$Whittaker model and $d_{\ppi}(S,\psi^D)=0$
by definition. The Theorem holds in this case.

 Now assume $\ppi^D=\ppi$. From Theorem~\ref{T:Wlocal}, for
$v\not\in S$, $\ppi_v=\Theta(\pi_v,\psi_v)$ is unramified, thus $\pi_v$
is unramified (see Proposition~4 in \cite{wa2}). As $\pi_v\otimes \chi_D$ 
is unramified for $v\not
\in S_1$, we can apply Theorem~\ref{T:21} to get:
\begin{equation}\label{E:proof4}
|d_{\pi\otimes\chi_D}(S_1,\psi)|^2
L^{S_1}(\pi\otimes\chi_D,1/2)=|d_{\ppi}(S_1,\psi^D)|^2.
\end{equation}
We will take a vector $\pp$ in the space of $\ppi$ so that $\pp=\otimes \pp_v$
with $\pp_v=\pp_{0,v}$ when $v\not\in S$. Using the equation (\ref{E:d2explicit}),
we get:
$$
d_{\ppi}(S_1,\psi^D)=\frac{|\W_{\pp}^D(e)|}{\|\pp\|}\prod_{v\in S_1}
\frac{\|\pp_v\|}{|\L^D_v(\pp_v)|}=d_{\ppi}(S,\psi^D)\prod_{v\in S_2}
\frac{\|\pp_{0,v}\|}{|\L^D_v(\pp_{0,v})|}.$$
That is
\begin{equation}\label{E:proof5}
|d_{\ppi}(S_1,\psi^D)|^2=|d_{\ppi}(S,\psi^D)|^2\prod_{v\in S_2}e(\pp_{0,v},\psi_v^D).
\end{equation}
\begin{lemma}\label{L:samed}
\begin{equation}\label{E:samed}
d_{\pi\otimes\chi_D}(S_1,\psi)=d_{\pi}(S_1,\psi).
\end{equation}
\end{lemma}
\begin{proof}
The space $V_{\pi\otimes\chi_D}$ consists of
the automorphic forms $\varphi(g)\chi_D(\det(g))$ where $\varphi\in V_{\pi}$.
Locally we will use the Whittaker model for $\pi_v$, that is,
$V_{\pi,v}$ consists of functions  $\varphi_v$ on $G$ satisfying 
$\varphi_v(n(x)g)=\psi_v(x)\varphi_v(g)$ with $G$ acting through right
translation. Then $L_v:\varphi_v\mapsto \varphi_v(e)$ is a Whittaker
functional. The space of $\varphi_v\chi_D(g)$ with 
$\varphi_v\in V_{\pi,v}$ is the Whittaker model of  $\pi_v\otimes \chi_D$,
with $L_v':\varphi_v\chi_D\mapsto \varphi_v(e)$ being a Whittaker functional.

Take $\varphi=\otimes_{v\in S_1}\varphi_v\otimes_{v\not\in S_1}\varphi_{0,v}$
in $V_{\pi}$ with $\varphi_{0,v}$ being unramified vector such that $L_v(\varphi_{0,v})=1$.
Then $\varphi\chi_D=\otimes_{v\in S_1}
\varphi_v\chi_D\otimes_{v\not\in S_1}\varphi_{0,v}\chi_D$
is in $V_{\pi\otimes \chi_D}$ with $\varphi_{0,v}\chi_D$ being unramified 
vector such that $L'_v(\varphi_{0,v}\chi_D)=1$. Clearly $L_v(\varphi_v)=L'_v(\varphi_v\chi_D)$
and $\|\varphi_v\|=\|\varphi_v\chi_D\|$ for all $v\in S_1$. From
(\ref{E:defW}) and (\ref{E:22}), we see $W_{\varphi}(e)=W_{\varphi\chi_D}(e)$
and
$\|\varphi\|=\|\varphi_v\|$. Thus from the explicit  formula (\ref{E:d1explicit}),
we get the identity (\ref{E:samed}).
\end{proof}

As $\pi_v$ is unramified for $v\not\in S$, as in (\ref{E:proof5}) we have
\begin{equation}\label{E:proof7}
|d_{\pi}(S_1,\psi)|^2=|d_{\pi}(S,\psi)|^2\prod_{v\in S_2}e(\varphi_{0,v},\psi_v).
\end{equation}
From equations (\ref{E:proof4}), (\ref{E:proof5}), (\ref{E:proof7})
and the Lemma, we get:
\begin{equation}\label{E:proof8}
|d_{\pi}(S,\psi)|^2L^{S_1}(\pi\otimes\chi_D,1/2)=|d_{\ppi}(S,\psi^D)|^2
\prod_{v\in S_2}\frac{e(\pp_{0,v},\psi_v^D)}
{e(\varphi_{0,v},\psi_v)}.
\end{equation}
For $v\in S_2$, $\pi_v$ is unramified and unitary, $\ppi_v=\Theta(\pi_v,\psi_v)$,
 the quotient
$\frac{e(\pp_{0,v},\psi_v^D)}
{e(\varphi_{0,v},\psi_v)}$ is given in Propositions \ref{P:81} and
\ref{P:83}; it equals $(|D|_vL(\pi_v\otimes\chi_D,1/2))^{-1}$. (This is
the only place we use the fact that $D$ is a square free integer).
As $S_1=S\cup S_2$,
we get from (\ref{E:proof8}):
$$|d_{\pi}(S,\psi)|^2L^S(\pi\otimes\chi_D,1/2)=|d_{\ppi}(S,\psi^D)|^2
\prod_{v\in S_2}|D|^{-1}_v.$$
Since $D$ is in $F^*$ and $|D|_v=1$ when $D\not\in S_1$,
$\prod_{v\in S_1}|D|_v=1$; thus $\prod_{v\in S_2}|D|^{-1}_v=
\prod_{v\in S}|D|_v$. We get (\ref{E:Wadele}).

 Now assume $D$ is such that $D\not\in F^{\epsilon_0}(\pi)$.
Then for some $v\in \Sigma$,
$(\frac{D_{v}}
{\pi_v})\not=\epsilon_{0,v}$; thus
$\Theta(\pi_v\otimes\chi_D,\psi_v^D)\not=\ppi_v$, and by Theorem~\ref{T:Wlocal},
$\ppi_v$ does not have a nontrivial $\psi_v^D$-Whittaker functional.
Therefore $\ppi$ does not have a nontrivial $\psi^D$-Whittaker model, and
$d_{\ppi}(S,\psi^D)=0$.
$\qed$

{\sc Proof of Theorem~\ref{T:KZadele}:} When $\prod_{v\in \Sigma}\epsilon_v
\not=\epsilon(\pi,1/2)$, we get for $D\in F^{\epsilon}(\pi)$,
$\epsilon(\pi\otimes\chi_D)=-1$, thus $L(\pi\otimes\chi_D,1/2)=0$.
In this case, $\ppi^{\epsilon}$ is zero dimensional, and by definition
$d_{\ppi^{\epsilon}}(S,\psi^D)=0$. When $\prod_{v\in \Sigma}\epsilon_v
=\epsilon(\pi,1/2)$, by
Theorem~\ref{T:Wglobal}, $\pi=S_{\psi}(\ppi^{\epsilon})$. Thus the Theorem
follows from equation (\ref{E:Wadele}).
$\qed$

{\sc Proof of Corollary~\ref{C:delta2}:} Let $\pi=S_{\psi}(\ppi)$
as in Theorem~\ref{T:Wglobal}. 
We note that $\epsilon(D,\pi)=\epsilon(D',\pi)$. Assume
$\ppi=\ppi^{\epsilon_0}$ for some $\epsilon_0$. When
$\epsilon(D,\pi)\not=\epsilon_0$, from  Theorem~\ref{T:Wadele},
$d_{\ppi}(S\cup S_{D,D'},\psi^D)=d_{\ppi}(S\cup S_{D,D'},\psi^{D'})=0$.
When $\epsilon(D,\pi)=\epsilon_0$, we follow the first part of the
proof of Theorem~\ref{T:Wadele}, replacing $S_1$ by $S\cup S_{D,D'}$ in
the argument. We get (see (\ref{E:proof4})):
$$|d_{\pi\otimes\chi_D}(S\cup S_{D,D'},\psi)|^2
L^{S\cup S_{D,D'}}(\pi\otimes\chi_D,1/2)=
|d_{\ppi}(S\cup S_{D,D'},\psi^D)|^2
,$$
$$|d_{\pi\otimes\chi_{D'}}(S\cup S_{D,D'},\psi)|^2
L^{S\cup S_{D,D'}}(\pi\otimes\chi_{D'},1/2)=
|d_{\ppi}(S\cup S_{D,D'},\psi^{D'})|^2.$$
The equation (\ref{E:delta2}) follows from the fact that 
$\chi_{D'}=\chi_D$ and Lemma~\ref{L:samed}. (The equation can
also be established directly as in the proof of Lemma~\ref{L:dpi}).
$\qed$

{\sc Proof of Theorem~\ref{T:conjecture}:} First note that for $D$ a
square free integer, for any fixed finite set of places $S$, the value of
$\prod_{v\in S-S_{\infty}}|D|_v$ lies in a finite set of positive numbers. Thus
we can disregard this quantity in our computation below.

Assume (\ref{E:rc2}) holds
for some $\alpha>0$. Given any $\pi$ irreducible cuspidal representation
of $PGL_2(\aa)$, we prove (\ref{E:lh}).
Let $D$ be a square free integer, by Theorem~\ref{T:KZadele},
there is a finite set $S$ of places, such that equation (\ref{E:KZadele})
holds for some $\ppi^{\epsilon}$. We may as well assume that $L(\pi\otimes
\chi_D,1/2)\not=0$. Then $d_{\pi^{\epsilon}}(S,\psi^D)\not=0$, and
we can find $\pp=\otimes_{v\in S}\pp_v\otimes_{v\not\in S}
\pp_{0,v}\in V_{\ppi^{\epsilon}}$ with
$\W^D_{\pp}(e)\not=0$. Recall
\begin{equation}\label{E:proof15}
|d_{\ppi^{\epsilon}}(S,\psi^D)|=|d_{\ppi^{\epsilon}}
(\pp,S_{\infty},\psi^D)|
\prod_{v\in S-S_{\infty}}\frac{\|\pp_v\|}{|\L^D_v(\pp_v)|}.
\end{equation}
For a given $v$, the value of $\frac{\|\pp_v\|}{|\L^D_v(\pp_v)|}$ does
depend on $D$ (as it depends on $\psi_v^D$). We put the dependence on $D$ in the notation and denote
the value as $\frac{\|\pp_v\|_D}{|\Lv^D(\pp_v)|}$. 
\begin{lemma}\label{L:finite}
For a fixed $v\in S-S_{\infty}$ and fixed $\pp_v$, there are only finitely many 
possible values of $\frac{\|\pp_v\|_D}{|\Lv^D(\pp_v)|}$.
\end{lemma}
\begin{proof}
As $|D|_v=1$ or $|D|_v=q^{-1}$, the value of $D$ lies in finitely many
cosets of $({\mathcal O_v}^*)^2$. Write $D=D_0\alpha^2$ with $\alpha\in
{\mathcal O_v}^*$, then we can let $\Lv^{D\delta_i}(\pp_v)=
\Lv^{D_0\delta_i}(\ppi_v(\v{\alpha})\pp_v)$
where $\delta_i$ are representatives of square classes of $F_v^*$.
Then from (\ref{E:26}), we get $\|\pp_v\|_D=\|\pp_v\|_{D_0}$. Meanwhile
$\Lv^D(\pp_v)=\Lv^{D_0}(\ppi_v(\v{\alpha})\pp_v)$. As $\ppi_v$ is admissible,
the set of $\{\ppi_v(\v{\alpha})\pp_v|\alpha\in {\mathcal O_v}^*\}$
is finite. There are only finitely many possible values of $\Lv^D(\pp_v)$ when
$D=D_0\alpha^2$. Thus ony finitely many possible values of the quotient
$\frac{\|\pp_v\|_D}{|\Lv^D(\pp_v)|}$.
\end{proof}
From the Lemma, there is a positive constant $c(\pp)$ depending
only on $\pp$, with 
$$\prod_{v\in S-S_{\infty}}\frac{\|\pp_v\|}{|\Lv^D(\pp_v)|}<c(\pp).$$
From  (\ref{E:proof15}) and (\ref{E:rc2}), we get:
$$|d_{\ppi^{\epsilon}}(S,\psi^D)|<<_{\ppi^{\epsilon},\pp,\alpha}
|D|_{S_{\infty}}^{\alpha-1/2}.
$$
Using (\ref{E:KZadele}), we get:
$$L(\pi\otimes\chi_D,1/2)<<_{\ppi^{\epsilon},\pp,\alpha}
|D|_{S_{\infty}}^{2\alpha}.$$
As the set of $\ppi^{\epsilon}$
is finite and determined by $\pi$, we get the inequality (\ref{E:lh}) for 
$\beta=2\alpha$.

Conversely, assume the inequality (\ref{E:lh}) holds for some
$\beta=2\alpha>0$, take any $\ppi\in \tilde{A_{00}}$ and 
$\pp\in V_{\ppi}$, we prove (\ref{E:rc2}). 
We may as well assume $\pp=\otimes_{v\in S} \pp_v\otimes_{v\not\in S}\pp_{0,v}$,
where $S$ is a large enough finite set of places. 
Let $\pi=S_{\psi}(\ppi)$.  From (\ref{E:sumKZ}) and our assumption, we get
$$|d_{\ppi}(S,\psi^D)|<<_{\pi,\alpha}|D|_{S_{\infty}}^{\alpha-1/2}.$$
Using the equation (\ref{E:proof15}) we get:
$$|d_{\ppi}(\pp,S_{\infty},\psi^D)|
\prod_{v\in S-S_{\infty}}\frac{\|\pp_v\|_D}{|\Lv^D(\pp_v)|}
<<_{\pi,\alpha}|D|_{S_{\infty}}^{\alpha-1/2}.$$
From Lemma~\ref{L:finite}, we see there is a constant $c'(\pp)>0$ such that
$$\prod_{v\in S-S_{\infty}}\frac{\|\pp_v\|_D}{|\Lv^D(\pp_v)|}>c'(\pp)$$
for all $D$. We get
$$|d_{\ppi}(\pp,S_{\infty},\psi^D)|
<<_{\pp,\pi,\alpha}|D|_{S_{\infty}}^{\alpha-1/2}.$$
As $\pi$ is determined by $\ppi$, the implied constant only depends on 
$\alpha$ and $\ppi$.
$\qed$

\section{Local factors: some examples}\label{S:local}

In this section, we compute the local factors $e(\varphi_v,\psi)$
and $e(\pp_v,\psi^D)$ in equation (\ref{E:main1})
 for some specific choices of the vectors
$\varphi_v$ and $\pp_v$.  The computation here is standard and fairly easy.
The result has
already been used in the proof of Theorem~\ref{T:Wadele}. It is
also used when we translate our formula into more explicit results
about cusp forms, (see the proof of Theorem~\ref{T:KZgeneral}).

In subsections \ref{S:principal}--\ref{S:special},
we assume $v$ is a nonarchimedean place, with odd residue characteristic.
For simplicity, we will assume $\psi_v$ 
has order 0, (and denote it simply by $\psi$),
 and  $D$ is either a unit or generates the prime ideal
in ${\mathcal O}_v$ 
at nonarchimedean places $v$.

The cases we consider are the following:

1. When $\pi_v$ is an unramified unitary representation of $G(F_v)$ where $v$
is an odd non-archimedean place. Then $\ppi^D_v=\Theta(\pi_v,\psi)$ for all $D\in F_v^*$ and
is an unramified unitary representation of $G'(F_v)$. We take
$\varphi_v$ and  $\pp_v$ to be the unramified vectors in $V_{\pi,v}$
and $V_{\ppi,v}$ respectively.

2. When $\pi_v$ is a holomorphic discrete series representation of $G({\Bbb R})$. Then
by Theorem~\ref{T:Wlocal}
$\ppi^D_v$ is either a holomorphic discrete series or an anti-holomorphic discrete
series representation of $G'({\Bbb R})$. We will only consider the case when $\ppi_v$
is the corresponding holomorphic discrete series. We take $\varphi_v$
and $\pp_v$ to be the minimal weight vectors in $V_{\pi,v}$
and $V_{\ppi,v}$ respectively.

3. When $\pi_v$ is a special representation of $G(F_v)$ where $v$
is a non-archimedean place. Then by Theorem~\ref{T:Wlocal} $\ppi^D_v$ could be either
$\ppi^+_v$ or $\ppi^-_v$. $\ppi^+_v$ is a special representation of $G'(F_v)$ while
$\ppi^-_v$ is a supercuspidal representation. We consider both cases. The vectors $\varphi_v$
and $\pp_v$ will be described in subsection~\ref{S:special}.

When we look at the cuspidal representations corresponding to the integral weight forms of
level $N$ and
half integral weight forms of level $4N$ with $N$ being
square free, the local components at infinite places and odd non-archimedean
places are of the form $\pi_v$
and $\ppi_v$ considered above. The situation at the even place is more subtle and will 
be considered in the next section.

\subsection{Some principal series at nonarchimedean places}\label{S:principal}

Let $\pi_v=\pi(\mu,\mu^{-1})$ be a unitary representation
with $\mu(x)=|x|_v^{s}$, $is\in {\Bbb R}$. Let $\ppi_v=\Theta(\pi_v,
\psi)$, it is $\ppi(\mu,\psi)$ by Proposition~4 of \cite{wa2}.
Note that $\ppi(\mu,\psi)$ is unramified.
We will let $\pp_v$ be  $\pp_{0,v}$ the unramified
vector in $V_{\ppi,v}$, and let $\varphi_v$ be $\varphi_{0,v}$
 the unramified vector in $V_{\pi,v}$.
Then $\pp_v$ and $\varphi_v$ are respectively $SL_2({\mathcal O}_v)$
 and $G({\mathcal O}_v)$ invariant functions in $\ppi(\mu,\psi)$ and
$\pi(\mu,\mu^{-1})$.

\begin{proposition}\label{P:81}
With above choices,
\begin{eqnarray*}
e(\varphi_v,\psi)&=&\frac{1+q^{-1}}{|1-q^{-2s-1}|^2},\\
e(\pp_v,\psi^D)&=&\left\{
\begin{array}{ll}
|D|_v^{-1}\frac{1+q^{-1}}{|1+q^{-1/2-s}\chi_D(\varpi)|^2}&\text{  when  }
|D|_v=1,\\
|D|_v^{-1}\frac{1+q^{-1}}{|1-q^{-2s-1}|^2}&\text{  when  }
|D|_v=q^{-1}.
\end{array}\right.
\end{eqnarray*}
Therefore  $\frac{e(\varphi_v,\psi)}{e(\pp_v,\psi^D)}=
|D|_vL(\pi_v\otimes\chi_D, 1/2)$.
\end{proposition}
\begin{proof}
We will use the Whittaker functional on $V_{\pi,v}$:
\begin{equation}\label{E:81}
L_v(\phi_v)=\int \phi_v(w n(x))\psi(-x)dx, \,\,\phi_v\in V_{\pi,v}
\end{equation}
The formula for spherical Whittaker function is well known, see \cite{cs}. The
formula for our case is available in \cite{g}.
We have
 $L_v(\pi_v(\u{a})\varphi_v)=0$ if
$|a|_v>1$; 
it equals 
$$q^{-m/2}\mu^{-1}(a)\frac{(1-q^{-2s-1})(1-q^{-2(m+1)s})}
{1-q^{-2s}}\varphi_v(e)$$
if $|a|_v=q^{-m}$ with $m\geq 0$. It is  easy to show
from (\ref{E:23}) that $\|\varphi_v\|^2$ equals
\begin{equation}\label{E:norm1}
|\frac{1-q^{-2s-1}}{1-q^{-2s}}|^2|\varphi_v(e)|^2\sum_{m=0}^{\infty}
q^{-m}(1-q^{-1})|1-q^{-2(m+1)s}|^2
=(1+q^{-1})|\varphi_v(e)|^2,
\end{equation}
 and thus the result
on $e(\varphi_v,\psi)$.

To compute $e(\pp_v,\psi^D)$, we use the $\psi^D-$Whittaker functional
\begin{equation}\label{E:84}
\Lv^D(\phi)=\int \phi(\tilde{w}\cdot \tilde{n}(x))\psi^D(-x)dx
\end{equation}
The formula for spherical Whittaker functions on the metaplectic
groups  is given in \cite{bfh}.
In our case it is an easy exercise to show
 when $|D|_v=1$, $\Lv^D(\pp_{0,v})=\pp_{0,v}(e)
(1+q^{-1/2-s}\chi_D(\varpi))$;
when $|D|_v=q^{-1}$, it equals $\pp_{0,v}(e)(1-q^{-1-2s})$.

The Hermitian form on $V_{\ppi,v}$ takes the form (\cite{bm} (9.19)):
\begin{equation}
(\phi,\phi')=\sum_{\delta_i\in F_v^*/(F^*_v)^2} 1/2\int\Lv^{D\delta_i}(\ppi(\v{a})\phi)
\overline{\Lv^{D\delta_i}(\ppi(\v{a})\phi')}|\delta_i|_v\frac{da}{|a|_v}.
\end{equation}
 The above form in turn equals: (\cite{bm} (9.18))
\begin{equation}\label{E:norm6}
\int |D|_v^{-1}\phi(\tilde{w}\cdot \tilde{n}(x))
\overline{\phi'(\tilde{w}\cdot \tilde{n}(x))}dx.\end{equation}
Use the second formula for Hermitian form to compute $\|\pp_{0,v}\|^2$.
From Iwasawa decomposition we get
\begin{equation}\label{E:norm3}
\|\pp_{0,v}\|^2=|D|_v^{-1}[\int_{|x|_v\leq 1}|\pp_{0,v}(e)|^2
dx+\int_{|x|_v>1}|x|_v^{-2}|\pp_{0,v}(e)|^2dx]
\end{equation}
which gives
$$\|\pp_{0,v}\|^2=|D|_v^{-1}(1+q^{-1})|\pp_{0,v}(e)|^2$$
This gives the result on $e(\pp_{0,v},\psi^D)$. The result on
the quotient $\frac{e(\varphi_v,\psi)}{e(\pp_{0,v},\psi^D)}$ follows
from the table on $L-$functions in \cite{g}. 
\end{proof}

\subsection{Complementary series at nonarchimedean places}\label{S:complementary}

Let  $\pi_v$  be  as in subsection~\ref{S:principal}, except that now
$\mu(x)=|x|_v^s\chi_{\tau}(x)$ with $\tau$ a unit in ${\mathcal O}_v$,
and $|s|<1/2, s\in {\Bbb R}$.
Let $\ppi_v=\Theta(\pi_v,\psi)$. 
Then as before $\ppi_v=\ppi(\mu,\psi)=\ppi(||_v^s,\psi^{\tau})$.
 We will choose
the vectors $\varphi_v$ and $\pp_v$ as in subsection~\ref{S:principal}.

\begin{proposition}\label{P:83}
With above choices,
\begin{eqnarray*}
e(\varphi_v,\psi)&=&\frac{1+q^{-1}}{(1-q^{-2s-1})(1-q^{2s-1})},\\
e(\pp_v,\psi^D)&=&
\left\{\begin{array}{ll}
|D|_v^{-1}\frac{1+q^{-1}}{(1+q^{-1/2-s}\chi_{D\tau}(\varpi))
(1+q^{s-1/2}\chi_{D\tau}(\varpi))}&
\text{  when  }|D|_v=1,\\ 
|D|_v^{-1}\frac{1+q^{-1}}{
(1-q^{2s-1})(1-q^{-1-2s})}&
\text{  when  }
|D|_v=q^{-1}.
\end{array}\right.
\end{eqnarray*}
 Therefore $\frac{e(\varphi_v,\psi)}{e(\pp_v,\psi^D)}=
|D|_vL(\pi_v\otimes\chi_D, 1/2)$.
\end{proposition}
\begin{proof} Retain the notations in the proof of Proposition~\ref{P:81}.
The formula for $L_v(\pi_v(\u{a})\varphi_v)$ still holds.
From (\ref{E:norm1}), one  gets $$\|\varphi_v\|^2=\frac{(1-q^{-1-2s})^2(1+q^{-1})}
{1-q^{2s-1}}|\varphi_v(e)|^2.$$ Thus we have the formula for
$e(\varphi_v,\psi)$.

The formula for $\Lv^D(\pp_{0,v})$ in the proof
of Proposition~\ref{P:81} remains valid.
The Hermitian form
however takes a more complicated form. If $z=\Delta^2 \delta$,
let
\begin{eqnarray*}
\lambda(z)&=&|\Delta|_v^{-2s-2}[(1-q^{-2s})^{-1}(1-q^{-1})+q^{s-1/2}
\chi_{\delta\tau}(\varpi)],
\text{  if  }|\delta|)v=1;\\
&=&|\Delta|_v^{-2s-2}q[(1-q^{-2s})^{-1}(1-q^{2s-1})],
\text{  if  }|\delta|_v=q^{-1}.
\end{eqnarray*}
Then $\lambda(z)=|z|_v^{s-1}\Delta(\psi,\tau,v)(z)$ where $\Delta(\psi,\tau,v)(z)$
is defined in \cite{bm} Proposition~9.8.
From equation (9.22) of \cite{bm}, the Hermitian form is:
\begin{equation}
(\phi,\phi')=\sum_{\delta_i\in F_v^*/(F_v^*)^2}1/2\int \Lv^{D\delta_i}(\ppi_v(\v{a})\phi)
\overline{\Lv^{D\delta_i}(\ppi_v(\v{a})\phi')}\frac{\lambda(z)}
{\lambda(D)}\frac{da}{|a|_v}. \end{equation}
From equation (9.21) of \cite{bm}, we see this form can also be written as:
\begin{equation}\label{E:89}
(\phi,\phi')=\lambda(D)^{-1}\int A\phi(\tilde{w}
\cdot\tilde{n}
(x))\overline{\phi'(\tilde{w}\cdot\tilde{n}(x))}dx.\end{equation}
where
$$A\phi(g)=\int\phi(\tilde{w}\cdot \tilde{n}(y) \cdot g)dy$$
Then $A\pp_{0,v}$ is the unique vector in the space of $\ppi(\mu^{-1}
\chi_{-1},\psi)$ fixed under $SL_2({\mathcal O}_v)$, with  $A\pp_{0,v}(e)=
\frac{1-q^{-1-2s}}{1-q^{-2s}}\pp_{0,v}(e)$. Using the Iwasawa
decomposition, we see
$$\|\pp_{0,v}\|^2=(1+q^{-1})
\lambda(D)^{-1}\frac{1-q^{-1-2s}}{1-q^{-2s}}|\pp_{0,v}(e)|^2.$$
This gives the formula for $e(\pp_{0,v},\psi^D)$.
The result on
the quotient $\frac{e(\varphi_v,\psi)}{e(\pp_{0,v},\psi^D)}$ follows
from the table on $L-$functions in \cite{g}. 
\end{proof}

\subsection{Special representations}\label{S:special}

\subsubsection{Description of $\ppi_v^+$ and $\ppi_v^-$}

 Let  $\mu_{\tau}(x)=|x|_v^{1/2}\chi_{\tau}(x)$, 
where $\tau$ is  in $F_v^*$. Let $\sigma^{\tau}=\sigma(\mu_{\tau}, \mu_{\tau}^{-1})$
be the special representation associated to the character $\mu_{\tau}$.
We will only consider the case when $|\tau|_v=1$.
The space of $\sigma^{\tau}$ is the
subspace of $\pi(\mu_{\tau},\mu_{\tau}^{-1})$ consisting of functions $\phi$
with
\begin{equation}\label{E:810}
\int\phi(wn(x))dx=0\end{equation}

From Theorem~\ref{T:Wlocal}, the set 
$\{\Theta(\sigma^{\tau}_v\otimes\chi_D,
\psi^D)\}$ consists of two elements. These two elements are
described in \cite{wa2}. When $D\tau$ is not a square,
$\Theta(\sigma^{\tau}_v\otimes\chi_D,
\psi^D)=\ppi^+$ is the special representation 
$\tilde{\sigma}^{\tau}(\psi)$. The space of this representation
is the subspace of $\ppi(\mu_{\tau},\psi)$ consisting of functions
$\phi$ satisfying
\begin{equation}\label{E:811}
\L_v^{\tau\Delta^2}(\phi)=\int \phi(\tilde{w}\cdot \tilde{n}(x))\psi(-\tau\Delta^2x)dx=0,
\text{  for all  }\Delta.\end{equation}
On the other hand, when $D\tau$ is a square, $\Theta(\sigma^{\tau}_v\otimes\chi_D,
\psi^D)=\ppi^-$ is a supercuspidal representation of $G'(F_v)$. More precisely,
it is the odd component of the Weil representation, denoted $r^-_{\psi^{\tau}}$.
The space of $r^-_{\psi^{\tau}}$ is the subspace of odd functions of $\cc(F_v)$,
with the action being: for $\Phi(z)\in \cc(F_v)$,
$\Phi(z)=-\Phi(-z)$,
\begin{equation}\label{E:weil1}
r^-_{\psi^{\tau}}(\tilde{n}(x))\Phi(z)=\psi^{\tau}(xz^2)\Phi(z),
\end{equation}
\begin{equation}\label{E:weil2}
r^-_{\psi^{\tau}}(\v{a})\Phi(z)=|a|_v^{1/2}\gamma(1,\psi^{\tau})/\gamma
(a,\psi^{\tau})\Phi(az),
\end{equation}
\begin{equation}\label{E:weil3}
r^-_{\psi^{\tau}}(\tilde{w})\Phi(z)=\gamma(1,\psi^{\tau})^2/\gamma(-1,\psi^{\tau})
\int \Phi(y)\psi^{\tau}(-2yz)dy.
\end{equation}

\subsubsection{The case of $\ppi^+_v$}

When $\tau D$ is not a square,  $\Theta(\sigma^{\tau}_v\otimes\chi_D,
\psi^D)$ is the special representation 
$\ppi^+_v=\tilde{\sigma}^{\tau}( \psi)$.

Define the Iwahori subgroup $K_0
\subset G({\mathcal O}_v)$ as the group
consisting of matrices $\mat{a}{b}{c}{d}$ with $|a|_v=|d|_v=1, |c|_v<1, |b|_v\leq 1$.
Recall in the introduction we defined an embedding of $SL_2({\mathcal O}_v)$ in 
$G'$ given by $g\mapsto (g,\kappa(g))$.
Let $K_0'$ be
the  image in $G'$
of the restriction of the splitting
to $K_0\cap SL_2$. 

Denote by $\text{char}(G({\mathcal O}_v))$ and Denote by $\text{char}(K_0)$
the characteristic functions if $G({\mathcal O}_v)$ and $K_0$ respectively.
Denote by $\text{char}(G'({\mathcal O}_v))$
a function on $G'$ with $\text{char}(G'({\mathcal O}_v))((g,\xi))$
equals 0 if $g\not\in SL_2({\mathcal O}_v)$, and
equals $\xi\kappa (g)$ otherwise.
Let $\text{char}(K'_0)(g,\xi)$ be the genuine
function on $G'$ that is 0 if $g
\not\in K_0$, and
equals $\text{char}(G'({\mathcal O}_v))(g,\xi)$ if $g\in K_0$.

\begin{lemma}\label{L:84}
Let $\varphi_v$ be a function in $\pi(\mu_{\tau },\mu_{\tau }^{-1})$ such that
it equals $
\text{char}(G({\mathcal O}_v))-(q+1)\text{char}(K_0)$ over
$G({\mathcal O}_v)$,
then $\varphi_v$ is in $V_{\sigma^{\tau },v}$ and  is fixed under $K_0$.

Let $\pp_v$ be a function in $\ppi(\mu_{\tau},\psi)$ such
that it equals $\text{char}(G'({\mathcal O}_v))-(q+1)\text{char}(K'_0)$ over
$G'({\mathcal O}_v)$,
then $\pp_v$ lies
in the space of $\tilde{\sigma}^{\tau}(\psi)$ and is fixed under $K'_0$.

The spaces of $K_0$ fixed vectors in $\sigma^{\tau }$ and
$K'_0$ fixed vectors in  $\tilde{\sigma}^{\tau}(\psi)$ are one dimensional.
\end{lemma}
\begin{proof}
We can consider the vectors in $\pi(\mu_{\tau },\mu_{\tau }^{-1})$
and $\ppi(\mu_{\tau},\psi)$ as functions on $G({\mathcal O}_v)$ and $SL_2({\mathcal O}_v)$
respectively.
Since $B\backslash G/K_0$ and $\tilde{B}\cap G'\backslash G'/K_0'$ both
have two elements,  the space of vectors in $\pi(\mu_{\tau },\mu_{\tau }^{-1})$
fixed by $K_0$ is two dimensional, with basis $\{\text{char}(G({\mathcal O}_v)),
\text{char}(K_0)\}$; the space of vectors in $\ppi(\mu_{\tau},\psi)$ fixed
by $K_0'$ is two dimensional, with basis $\{\text{char}(G'({\mathcal O}_v)),
\text{char}(K_0')\}$. 

When $\phi_1$ is the vector corresponding to $\text{char}(G({\mathcal O}_v))$,
 the integral $\int \phi_1(wn(x))dx$ equals $1+q^{-1}$.
When $\phi_2$ is the vector corresponding to $\text{char}(G(K_0))$,
the integral $\int \phi_2(wn(x))dx$ equals $q^{-1}$.
Thus $\varphi_v$ satisfies the condition (\ref{E:810}) and
generates the one dimensional space fixed by $K_0$ in $\sigma^{\tau }$.

Next we  compute $\Lv^z(\pp_v)$ with $z\in F_v^*$.
Let $\phi_1'$ be the vector corresponding to $\text{char}(G'({\mathcal O}_v))$,
then from (\ref{E:84}) and the Iwasawa decomposition, we get
$$\Lv^z(\phi_1')=\int_{|x|_v\leq 1}\psi(-zx)dx+\sum_{r=1}^{\infty}
\int_{|x|_v=q^r}\tilde{\gamma}(x,\psi)|x|_v^{-3/2}\chi_{\tau}(x)[-1,x]\psi(-zx)dx.$$
When $|x|_v=q^r$ with $r$  even, with our assumption on $\tau$ being
a unit, $$\tilde{\gamma}(x,\psi)\chi_{\tau}(x)[-1,x]=1.$$ When 
$|x|_v=q^r$ with $r$ odd, 
$$\tilde{\gamma}(x,\psi)\chi_{\tau}(x)[-1,x]
=\gamma(x,\psi)[\tau,\varpi].$$
It is a simple calculation to get the following result: write $z=\delta\Delta^2$ with
$|\delta|_v=1$ or $q^{-1}$,
$$
\Lv^{\delta\Delta^2}(\phi_1')=\left\{
\begin{array}{ll}
0&\text{when } |\Delta|_v>1,\\
1+q^{-1}+|\Delta|_v(q^{-1}[\tau\delta,\varpi]-q^{-1})&
\text{when }|\delta|_v=1,|\Delta|_v\leq  1,\\
q^{-1}-|\Delta|_v(q^{-1}+q^{-2})&
\text{when }|\delta|_v=q^{-1},|\Delta|_v\leq  1.
\end{array}\right.$$

Let $\phi_2'$ be the vector corresponding to $\text{char}(K_0')$. Then
$$\Lv^z(\phi_1'-\phi_2')=\int_{|x|_v\leq 1}\psi(-zx)dx.$$ 
We get
$$
\Lv^{\delta\Delta^2}(\phi_1')=\left\{
\begin{array}{ll}
0&\text{when } |\Delta|_v>1,\\
q^{-1}+|\Delta|_v(q^{-1}[\tau\delta,\varpi]-q^{-1})&
\text{when }|\delta|_v=1,|\Delta|_v\leq  1,\\
q^{-1}-|\Delta|_v(q^{-1}+q^{-2})&
\text{when }|\delta|_v=q^{-1},|\Delta|_v\leq  1.
\end{array}\right.$$

The formula for $\Lv^{z}(\pp_v)$ is
\begin{equation}\label{E:specialwhittaker}
\Lv^{\delta\Delta^2}(\pp_v)=\left\{
\begin{array}{ll}
0&\text{when } |\Delta|_v>1,\\
2|\Delta|_v&
\text{when }|\delta|_v=1,|\Delta|_v\leq  1, \delta\tau \text{ is not a square},\\
0& \text{when }|\delta|_v=1,|\Delta|_v\leq  1, \delta\tau \text{ is a square},\\
|\Delta|_v(q^{-1}+1)&
\text{when }|\delta|_v=q^{-1},|\Delta|_v\leq  1.
\end{array}\right.
\end{equation}
It is now clear that $\Lv^z(\pp_v)$ satisfies the condition
(\ref{E:811}).  The vector $\pp_v$ generates the space of $K_0'$
fixed vectors in $\tilde{\sigma}^{\tau}(\psi)$ which from the above
formulas is clearly one dimensional.
\end{proof}

\begin{proposition}\label{P:85}
Assume $D\tau$ is not a square. Let $\varphi_v$ and $\pp_v$
be the vectors in Lemma~\ref{L:84}.
Then
\begin{eqnarray*}
e(\varphi_v,\psi)&=&\frac{1}{1+q^{-1}},\\
e(\pp_v,\psi^D)&=&
\left\{\begin{array}{ll}
1/2& \text{  when }|D|_v=1,\\
q/(1+q^{-1})& \text{  when }|D|_v=q^{-1}.
\end{array}\right.
\end{eqnarray*}
Therefore
$$\frac{e(\varphi_v,\psi)}{e(\pp_v,\psi^D)}=\left\{\begin{array}{ll}
2L(\pi_v\otimes\chi_D, 1/2)|D|_v & \text{  when }|D|_v=1,\\
L(\pi_v\otimes\chi_D, 1/2)|D|_v& \text{  when }|D|_v=q^{-1}.
\end{array}\right.$$
\end{proposition}
\begin{proof} We can use the Iwasawa decomposition
to compute $L_v(\sigma^{\tau }(\u{a})\varphi_v)$
where $L_v$ is defined as in (\ref{E:81}).
 We will skip the details. One gets
$$ L_v(\sigma^{\tau }(\u{a})\varphi_v)=\left\{
\begin{array}{ll}
0&\text{when }|a|_v>1,\\
(1+q^{-1})|a|_v\chi_{\tau }(a)&\text{when }|a|_v\leq 1.
\end{array}\right.
$$
 Thus
$\|\varphi_v\|^2=(1+q^{-1})$ from (\ref{E:23}) and we get
the value of $e(\varphi_v,\psi)$.

Assume $|D|_v=1$. From (\ref{E:specialwhittaker}) $\Lv^D(\pp_v)=2$.
To find $\|\pp_v\|$,  we use the Hermitian form (9.23) in
\cite{bm}:
\begin{equation}\label{E:Hspecial}
(\phi,\phi')=\sum_{b=D,D\varpi,\tau\varpi}c_b/2\int \Lv^b(\tilde{\sigma}(\v{a})\phi)
\overline{\Lv^b(\tilde{\sigma}(\v{a})\phi')}\frac{da}{|a|}_v
\end{equation}
where $\tilde{\sigma}=\tilde{\sigma}^{\tau}(\psi)$ and $c_b=1$
when $b=D$, and $c_b=2q^{-1}/(1+q^{-1})$ when $b=D\varpi, \tau\varpi$. Using
this formula, the formula (\ref{E:specialwhittaker})
for $\Lv^z(\pp)$ 
 and the fact that $$|\Lv^b(\tilde{\sigma}(\v{a})\pp)|=|a|_v^{1/2}
|\Lv^{a^2b}(\pp)|,$$
we get 
$$\|\pp\|^2=\int_{|a|_v\leq 1} 2|a|_v^3\frac{da}{|a|_v}+\frac{2q^{-1}}{1+q^{-1}}
\int_{|a|_v\leq 1} (1+q^{-1})^2|a|_v^3\frac{da}{|a|_v}=2.$$
Thus we get the formula for
 $e(\pp_v,\psi^D)$ when $|D|_v=1$

When $|D|_v=q^{-1}$,. from (\ref{E:specialwhittaker}) $\Lv^{D}(\pp_v)=1+q^{-1}$.
The computation of $\|\pp_v\|$ goes (\ref{E:26}) the Hermitian form changes
from (\ref{E:Hspecial}) to 
$$(\phi,\phi')=\sum_{b=D,D\tau,\delta}c'_b/2\int \Lv^b(\tilde{\sigma}(\v{a})\phi)
\overline{\Lv^b(\tilde{\sigma}(\v{a})\phi')}\frac{da}{|a|_v}
$$
where $\delta$ is a unit in ${\mathcal O}_v$ such that $\delta\tau$ is not a square. Here
$c'_D=c'_{D\tau}=1$ and $c'_{\delta}=\frac{1+q^{-1}}{2q^{-1}}$. 
We get $\|\pp_v\|=(1+q^{-1})/q^{-1}$ and the formula for $e(\pp_v,\psi^D)$.

The claim on the quotient follows from the following formulas for $L-$values (see \cite{g}).
When $|D|_v=1$ with $\tau D$ not a unit,
$L(\pi_v\otimes \chi_D,1/2)=L(\sigma^{\tau D},1/2)=(1+q^{-1})^{-1}.$
When $|D|_v=q^{-1}$, $L(\pi_v\otimes \chi_D,1/2)=L(\sigma^{\tau D},1/2)=1$.
\end{proof}

\subsubsection{The case of $\ppi_v^-$}

When $\tau D$ is a square,  $\Theta(\sigma^{\tau}_v\otimes\chi_D,
\psi^D)$ is the supercuspidal representation $\ppi_v^-=r^-_{\psi^{\tau}}$.
Recall by our assumption $\tau$ and $D$ are both units in ${\mathcal O}_v$.
We will let $\varphi_v$ to be the vector defined in Lemma~\ref{L:84}. Next we describe
a vector $\pp_v$ in the space of $r^-_{\psi^{\tau}}$.

Let $K_{00}$ be the subgroup of $K_0$ consisting 
of matrices $\mat{a}{b}{c}{d}$ with $|a|_v=|d|_v=1, |c|_v<q^{-1}, |b|_v\leq 1$.
Let $K_{00}'$ be the image of $K_{00}\cap SL_2$ embedded in $G'$. Then 
$K_{00}'=\{(\sigma,\kappa(\sigma))|\sigma\in K_{00}\}$. Let
$\chi$ be any odd character of ${\mathcal O}_v^*$ that is
trivial on $1+P$, $(\chi(-1)=-1)$, then $\chi$ defines a character on
$K_{00}'$ by:
$$\chi(\sigma,\kappa(\sigma))\mapsto \chi(d), \,\sigma=\mat{a}{b}{c}{d}\in K_{00}.$$

Let $char(X)$ denote the characteristic function of a subset $X$ in
$F_v$.
\begin{proposition}\label{L:weil1}
The space of vectors $\Phi$ in $r^-_{\psi^{\tau}}$ satisfying 
$r^-_{\psi^{\tau}}(k)\Phi=\chi(k)\Phi$, $k\in K_{00}'$ is 
one dimenisonal. It is generated by the element
$$\Phi_{\chi}(z)=\sum_{b\in {\mathcal O}_v^*/1+P}\chi^{-1}(b)\text{char}(b+P)(z).$$
\end{proposition}
\begin{proof}
Let $\Phi$ be a vector satisfying the relation in the Proposition.
With our assumptions on $\psi,\tau$ and the place $v$, the equation
(\ref{E:weil2}) gives
$$r^-_{\psi^{\tau}}(\v{z})\Phi(1)=\Phi(z)=\chi(z^{-1})\Phi(1), \,\,z\in {\mathcal O}_v^*.$$
Thus $\Phi_0=\Phi-\Phi(1)\Phi_{\chi}$ vanishes over the set ${\mathcal O}_v^*$.
Next if $z\not \in {\mathcal O}_v$, from (\ref{E:weil1}), for $x\in
{\mathcal O}_v$:
$$r^-_{\psi^{\tau}}(\tilde{n}(x))\Phi(z)=\psi(\tau zx^2)\Phi(z)=\Phi(z),$$
thus $\Phi(z)=0$. We get $\Phi_0$ is supported on $P$.

We now show $\Phi_{\chi}$ satisfies the relation in the Proposition.
\begin{lemma}\label{L:weil2}
For $z\in F_v$, $x\in {\mathcal O}_v$,
\begin{equation}\label{E:weil5}
r^-_{\psi^{\tau}}(\mat{1}{}{x\varpi^2}{1},1)\text{char}(z+P)
=\text{char}(z+P).
\end{equation}
\end{lemma}
\begin{proof}
We use the fact that
$$(\mat{1}{}{x\varpi^2}{1},1)=\tilde{w}\cdot (-e,1)\cdot \tilde{n}(-x\varpi^2)
\cdot \tilde{w}.$$
From (\ref{E:weil1}), (\ref{E:weil2}) and
(\ref{E:weil3}) the left hand side of (\ref{E:weil5}) is
$$r^-_{\psi^{\tau}}(\mat{1}{}{x\varpi^2}{1},1)\text{char}(z+P)(a)=
\int\int  \text{char}(z+P)(u)\psi^{\tau}(xy^2\varpi^2+2uy-2ay)dudy.$$
The integral over $u$ is nonzero only when $y\in P^{-1}$, in which case
$\psi^{\tau}(xy^2\varpi^2)=1$ and from
the Fourier inversion formula, the above integral just equals $\text{char}(z+P)(a)$.
\end{proof}

 It
is easy to check using (\ref{E:weil1}) and (\ref{E:weil2})
that $r^-_{\psi^{\tau}}(k))\Phi_{\chi}=\chi(k)\Phi_{\chi}
$ when $k\in K_{00}'\cap \tilde{B}$. From
the above Lemma we get for
$x\in {\mathcal O}_v$, 
$$r^-_{\psi^{\tau}}(\mat{1}{}{x\varpi^2}{1},1)\Phi_{\chi}=\Phi_{\chi}.$$
Since these group $K_{00}'$ is generated by the elements in 
$K_{00}'\cap \tilde{B}$ and $\{(\mat{1}{}{x\varpi^2}{1},1)|x\in {\mathcal O}_v\}$,
we see $\Phi_{\chi}$ satisfies the relation in the Proposition.

Thus $\Phi_0=\Phi-\Phi(1)\Phi_{\chi}$ satisfies the relation in the Proposition
and is supported over $P$. 
To finish the proof, we need to show such a function is identically 0.
From the proof of Lemma~\ref{L:weil2}, we get:
$$\int_{x\in {\mathcal O}_v}r^-_{\psi^{\tau}}(\mat{1}{}{x\varpi^2}{1},1)\Phi_0(a)
=\Phi_0(a)$$
$$=\int_{x\in {\mathcal O}_v}\int\int
 \Phi_0(u)\psi^{\tau}(xy^2\varpi^2+2uy-2ay)dudydx.$$
For the integration over $x$ to be nonzero, $y\in P^{-1}$, in which case
$\Phi_0(u)\psi^{\tau}(2uy)=\Phi_0(u)$. Thus
the above integral equals
$$
\int_{x\in {\mathcal O}_v}\int_{y\in P^{-1}}\int
 \Phi_0(u)\psi^{\tau}(-2ay)dudydx$$
which equals a constant times $\text{char}(P)$. Thus for $a\in P$,
$\Phi_0(a)=\Phi_0(0)$.
Since $\Phi_0$ is an odd function, $\Phi_0$ vanishes over $P$, thus
vanishes identically.
\end{proof}

The representation $r^-_{\psi^{\tau}}$ is a {\em distinguished}
representation, in the sense that it only has nontrivial Whittaker
functionals for $\psi^{\delta}$ with $\delta$ in the same square
class as $\tau$. Assume $D=\tau \alpha^2$, 
we can define $\Lv^{D}$ by setting
$$\Lv^{D}(\Phi)=\Phi(\alpha).$$
Then the Hermitian form is just:
$$(\Phi_1,\Phi_2)=1/2\int \Lv^{D}(r^-_{\psi^{\tau}}(\v{a})\Phi_1)
\overline{\Lv^{D}(r^-_{\psi^{\tau}}(\v{a})\Phi_1)}\frac{da}{|a|}_v.$$
Clearly $\Lv^D(\Phi_{\chi})=\chi(\alpha^{-1})$ and by (\ref{E:weil2}),
$$|\Lv^D(r^-_{\psi^{D}}(\v{a})\Phi_{\chi})|=|a|_v^{1/2}|\Phi_{\chi}(a\alpha)|$$
which equals $1$ when $|a|_v=1$ and $0$ otherwise. Thus
$$\|\Phi_{\chi}\|^2=1/2\int_{|a|_v=1} \frac{da}{|a|_v}=(1-q^{-1})/2.$$
Thus $e(\Phi_{\chi},\psi^D)=(1-q^{-1})/2$. Note that from \cite{g}
$$L(\pi_v\otimes\chi_D,1/2)=L(\sigma^{\tau D},1/2)=L(\sigma^{1},1/2)=(1-q^{-1})^{-1}.$$
 We have
\begin{proposition}\label{P:weil1}
Let $\varphi_v$ be a vector in $\pi_v=\sigma^{\tau }$ given by Lemma~\ref{L:84}.
Let $\pp_v$ be $\Phi_{\chi}$ as in Proposition~\ref{L:weil1}. Then
$e(\varphi_v,\psi)=(1+q^{-1})^{-1}$, $ e(\pp_v,\psi^D)=(1-q^{-1})/2$ and
$\frac{e(\varphi_v,\psi)}{ e(\pp_v,\psi^D)}=2(1+q^{-1})^{-1}L(\pi_v\otimes\chi_D,1/2)$.
\end{proposition}

\subsection{Holomorphic discrete series}\label{S:holomorphic}

Let $F_v={\Bbb R}$. Let $\pi_v$ be the discrete series $\sigma(\mu,\mu^{-1})$ 
(\cite{wa2}) where
$\mu(x)=|x|^{s/2}(\text{sgn }x)^{(s+1)/2}$, $k=\frac{s+1}{2}$ being
a positive integer.  Then according to  Theorem~\ref{T:Wlocal}, $\Theta(\pi\otimes
\chi_D,\psi^D)$ can be one of the following two representations:
$\ppi_v=\ppi_v^+=\Theta(\pi_v,\psi)$ and $\ppi_v^-=\Theta(\pi_v\otimes \text{ sgn},\psi^{-1})$. 
We note in this case $\pi_v\cong \pi_v\otimes \text{sgn }$.

We now assume $\psi(x)=e^{2\pi inx}$ with $n$ a positive
integer. Then   $\ppi_v$ is a holomorphic discrete series 
$\tilde{\sigma}(\mu)$  while $\ppi_v^-=\tilde{\sigma}(\mu\text{ sgn})$ is 
an antiholomorphic discrete series, (\cite{wa2}).

Let $\varphi_v$ and $\pp_v$ be a vector of minimal weight in
$V_{\pi,v}$ and $V_{\ppi,v}$ respectively. These vectors
are determined up to a scalar. We have:
\begin{proposition}\label{P:86}
Let $D$ be a positive integer. With the above notations:
\begin{eqnarray*}
e(\varphi_v,\psi)&=&e^{4\pi n}(4\pi n)^{-2k}\Gamma(2k)\\
e(\pp_v,\psi^D)&=& 2e^{4\pi nD}(4\pi nD)^{-(1/2+k)}\Gamma(1/2+k)
\end{eqnarray*}
Therefore $\frac{e(\varphi_v,\psi)}{e(\pp_v,\psi^{D})}=\frac{1}{2}
e^{4\pi n(1-D)}
D^{1/2+k}n^{1/2-k}\pi^{-k}(k-1)!$.
\end{proposition}
\begin{proof}
 From \cite{g}, we see the Whittaker model for $\varphi_v$
has the form:
$$L_v(\pi_v(\u{a}){\varphi_v})=\left\{\begin{array}{cc}\alpha a^{k}e^{-2\pi na}&
a>0\\0&a<0\end{array}\right.$$
where $\alpha$ is some nonzero constant which we may as well fix
to be 1.
With this model, we  get from (\ref{E:23})
$$\|\varphi_v\|^2=\int_{a>0}a^{2k}e^{-4\pi
na}d^{*}a$$
which equals
$(4\pi n)^{-(2k)}\Gamma(2k)$. This gives the result for $e(\varphi_v,
\psi)$.

From \cite{wa} p.24, we see the Whittaker model for $\pp_v$ with
respect to $\psi^D$ has the form:
$$\Lv^D(\ppi_v(\v{a})\pp_v)=\alpha \omega(sgn(a))|a|^{1/2+k}e^{-2\pi nDa^2}$$
where $\omega$ is the central character of $\tilde{\sigma}(\mu)$.
We will again let $\alpha=1$. Since in this case, $\ppi_v$ is
{\em distinguished}, i.e., the Whittaker functional for $\psi^{z}$
is always trivial when $z<0$, the Hermitian form in (\ref{E:26}) simplifies to:
$$(\phi,\phi')=\int \Lv^D(\ppi_v(\v{a})\phi)\overline{\Lv^D(\ppi_v(\v{a})
\phi')}d^{*}a$$
Apply the formula to compute $\|\pp_v\|$, we find that
$$\|\pp_v\|^2=2(4\pi nD)^{-(1/2+k)}\Gamma(1/2+k).$$ The result for
$e(\pp_v,\psi^D)$ follows. The assertion  on the quotient follows  from
the formula ($k$ an integer):
$$\Gamma(2k)=\pi^{-1/2}2^{2k-1}\Gamma(k)\Gamma(k+1/2), \, \Gamma(k)=(k-1)!.$$
\end{proof}

\section{Cusp forms over $\Bbb Q$}\label{S:dictionary}

We will apply the results in \S~\ref{S:mainresult} to the case of holomorphic cusp forms over ${\Bbb Q}$.
Fix the additive character $\psi$ as follows:  if
$x\in {\Bbb R}$,  $\psi(x)=e^{2\pi ix}$; at a rational prime $p$,
if $x\in {\Bbb Q}_p$, choose $\hat{x}\in  {\Bbb Q}$ so that $|x-\hat{x}|_p\leq 1$, and
set $\psi(x)=e^{-2\pi i\hat{x}}$. Denote by $\tilde{\gamma}(x)$ the number
$\tilde{\gamma}(x,\psi)$. We denote by $|D|_v$ the metric at a place
$v$, and $|D|$ the absolute value of $D$ which equals $|D|_{\infty}$.

We first recall the correspondence between the cusp forms and automorphic
representations. Our reference is \cite{wa1} section III. The main result
in this section is the choise of  a one dimensional subspace in a
two dimensional subspace of $V_{\ppi_2}$. This choice
is closely related to the definition of  the Kohnen space of
half-integral weight forms.

\subsection{The dictionary: integral weight form}\label{S:integer}

 Let $\Gamma_0(N)=\{
\mat{a}{b}{c}{d}\in SL_2({\Bbb Z})| c\equiv 0 (N)\}$. Let
$S_{2k}(N)$ be the space of cusp form of weight $2k$ on $\Gamma_0(N)$ (of
 level $N$), and with trivial
character. Assume from now on that $N$ is odd and
square free.  Let $f\in S_{2k}(N)$ be a newform.
Then $f$ determines a vector in the space of automorphic forms on
$GL_2(\aa_{\Bbb Q})$ by
$f\mapsto \varphi=s(f)$. The map $s(f)$ is defined as
follows. For $g_{\infty}=\mat{a}
{b}{c}{d}\in GL_2(\Bbb R)$, let
$$f|_{g_{\infty}}(z)=f(\frac{az+b}{cz+d})
(cz+d)^{-2k}$$
Consider $g_{\infty}$ as an element
$(g_{\infty},e,e,\ldots)$ in
$GL_2(\aa_{\Bbb Q})$, then $\varphi(g_{\infty})=f|_{g_{\infty}}
(i)$, and $\varphi(\gamma gk)=\varphi(g)$ whenever
$\gamma\in GL_2({\Bbb Q})Z(\aa_{\Bbb Q})$, and $k\in
\prod_{p\not|N}GL_2(\oo_p)\prod_{p|N}K_{0,p}$.

Then $\varphi$ is a vector in the space of 
an irreducible  cuspidal representation $\pi$ of
$GL_2(\aa_{\Bbb Q})$, with trivial central character. The representation
$\pi=\otimes\pi_v$, and $\varphi=\otimes_v\varphi_v$ can be
described as follows:

(\ref{S:integer}.1). When $v=\infty$, $\pi_v$ is the discrete series
$\sigma(\mu_{\infty},\mu_{\infty}^{-1})$ as in subsection~\ref{S:holomorphic}, 
with $\mu_{\infty}(x)=|x|_v^{k-1/2}(\text{sgn }x)^k$. The vector
$\varphi_{\infty}$ is a minimal weight vector.

(\ref{S:integer}.2). When $v$ is $p-$adic, $p$ not dividing $N$, then
$\pi_v=\pi(\mu_v, \mu_v^{-1})$ with $\mu_v$ an unramified
character, and $\varphi_v$ is an unramified
 vector.

(\ref{S:integer}.3). When $v$ is $p-$adic, $p|N$, then $\pi_v$ is a special
representation $\sigma^{\tau_v}$ as in subsection~\ref{S:special}, where
$\tau_v$ is a unit in ${\Bbb Z}_p$. 
 Then $\varphi_v$ is the vector described in Lemma~\ref{L:84}.

Conversely, given an irreducible cuspidal automorphic
representation $\pi=\otimes\pi_v$ with local components as described
in (\ref{S:integer}.1)--(\ref{S:integer}.3), pick $\varphi$ as above (which is unique up to
scalar multiple), then $\varphi$ is a scalar multiple of $s(f)$ for some newform $f$ 
in $S_{2k}(N)$.

If $f(z)=\sum_{n=1}^{\infty}a(n)e^{2\pi i nz}$, then
 $a(1)=e^{2\pi}W_{s(f)}(e)$. As we assume $a(1)=1$, for
$\varphi=s(f)$, $W_{\varphi}(e)=e^{-2\pi}$. 

\subsection{Dictionary: half integral weight form.}\label{S:half}

Assume now that $k$ is a nonnegative integer. Let $N$ be a positive odd
integer. Let $\chi$ be a Dirichlet character mod $4N$ such that
$\chi(-1)=1$. Assume $4N=\prod_{p|4N}p^{v(p)}$, then
 $({\Bbb Z}/4N)^*\cong \prod_{p|N}({\Bbb Z}/p^{v(p)})^*$, and $\chi$
can be decomposed into a product of characters $\chi_{(p)}$
of $({\Bbb Z}/p^{v(p)})^*$ under this isomorphism. We can trivially extend
$\chi_{(p)}$ to a character of ${\Bbb Z}_p^*$. 

Let $S'_{k+1/2}(4N,\chi)$ be the space of holomorphic
cusp forms of  weight $k+1/2$, level $4N$ and character $\chi$. The functions in
the space satisfies: \cite{wa1}
\begin{equation}
g(\frac{az+b}{cz+d})=j(\sigma,z)^{2k+1}\chi(d)g(z),\,
 \sigma=\mat{a}{b}{c}{d}\in SL_2({\Bbb Z}), \, 4N|c.
\end{equation} 
Here $$j(\sigma,z)=\theta(\frac{az+b}{cz+d})/\theta(z), \,\,\theta(z)=
\sum_{n={-\infty}}^{\infty}e^{2\pi in^2z}.$$ 

Let $\tilde{A}'_{k+1/2}(4N,\chi)$ be the space generated by vectors 
$\pp=\otimes_v\pp_v$ in the space of cuspidal automorphic forms on $\tilde{SL}_2(\aa)$
satisfying:

(\ref{S:half}.1) When $v=\infty$, $\pp_v$ is a minimal weight vector in
the space of a holomorphic discrete series $\tilde{\sigma}(\mu_{\infty})$
where $\mu_{\infty}(x)=|x|^{k-1/2}(\text{sgn }x)^{k}$. 

(\ref{S:half}.2) When $v$ is $p-$adic, $p$ not dividing $2N$, then
$\pp_p$ is the unramified vector in the space of $\ppi(\mu_v,\psi)$ where
$\mu$ is an unramified character.

(\ref{S:half}.3) When $v$ is $p-$adic, $p\not=2$, and $p|N$, then
$\pp_p$ is a vector in the space of some $\ppi_p$ such that $\ppi_p(\sigma,\kappa(\sigma))\pp_p
=\chi_{(p)}(d)\pp_v$ whenever $\sigma=\mat{a}{b}{c}{d}\in SL_2({\Bbb Z}_p)$
with $|c|_p\leq |N|_p$.

(\ref{S:half}.4) When $v$ is $p-$adic with $p=2$, $\pp_2$ is a vector
in the space of some $\ppi_2$ such that $\ppi_2(\sigma,1)\pp_2
=\tilde{\epsilon}_2(\sigma)\chi_{(2)}\chi_{-1}^{k}(d)\pp_2$ whenever
$\sigma=\mat{a}{b}{c}{d}\in SL_2({\Bbb Z}_2)$
with $|c|_2\leq |4N|_2$. Here
$$\tilde{\epsilon}_2(\sigma)=\left\{\begin{array}{ll}
\tilde{\gamma}(d)[c,d]& \text{when  }c\not=0\\
\tilde{\gamma}(d)^{-1}&\text{when  }c=0.
\end{array}
\right.$$

The Proposition~3 of \cite{wa1} establishes a bijection 
 from $S'_{k+1/2}(4N,\chi)$ to $\tilde{A}'_{k+1/2}(4N,\chi)$. 
The bijection is given by $g(z)\mapsto
\pp=t(g)$, where $t(g)$ is the unique function on $\tilde{SL}_2(\aa)$
that is continuous and left invariant under $SL_2({\Bbb Q})$ and satisfies: 
$$t(g)(\mat{\sqrt{y}}{x/\sqrt{y}}{0}{1/\sqrt{y}}\mat{\cos\theta}{\sin\theta}
{-\sin\theta}{\cos\theta},1)= y^{k/2+1/4}e^{i(k+1/2)\theta}g(x+yi),$$
where $y>0,x\in {\Bbb R}$ and $-\pi< \theta\leq \pi$.

The relation between the Whittaker functionals of $t(g)$ and
the Fourier coefficients of $g(z)$ is given by follows:
 From Lemma~3 of \cite{wa1},
we get when $g(z)=\sum_{n=1}^{\infty}c(n)e^{2\pi i nz}$, 
\begin{equation}\label{E:fc2}
c(n)=e^{2\pi n}\W^n_{t(g)}(e).
\end{equation}

{\bf Remark on Petersson norm}: If $f$ is a cusp form of weight $k\in \frac{1}{2}{\Bbb Z}$
on ome subgroup $\Gamma$ of finite index in $\Gamma_1=SL_2({\Bbb Z})$, we define 
as usual the norm of $f$ to be 
$$<f,f>=\frac{1}{[\Gamma(1):\Gamma]}\int_{\Gamma\backslash {\mathcal H}}|f(z)|^2y^{k-2}dxdy$$
where $z=x+iy$ and ${\mathcal H}$ is the upper half plane. Then
\begin{lemma}\label{L:petersson}
For $\varphi=s(f)$ and $\pp=t(g)$ as above:
\begin{equation}\label{E:petersson}
\frac{\|\varphi\|^2}{<f,f>}=\frac{|\pp\|^2}{<g,g>}.
\end{equation}
\end{lemma}
\begin{proof}
It is well known that $\|\varphi\|^2=<f,f>$ and $\|\pp\|^2=<g,g>$ when we use the
following Haar measures $d'$ on $GL_2$ and $SL_2$ instead of the one given in the
introduction. When $v$ is a nonarchimedean place, choose the measure $d'$ on $GL_2$ so that
$G(\oo_v)$ has volume 1;  choose the measure $d'$ on $SL_2$
so that $SL_2(\oo_v)$ has volume 1. When $v$ is the infinite place, 
let $k(\theta)=\mat{\cos\theta}
{\sin\theta}{-\sin\theta}{\cos\theta}$. 
From the Iwasawa decomposition any
$g\in GL_2^+({\Bbb R})$ (the subgroup with positive determinant)
 can be written uniquely as $g=z(c)n(x)\mat{a}{}{}{a^{-1}}k(\theta)$
with $c\in {\Bbb R}^{*}, x\in {\Bbb R}, a>0$ and $0\leq \theta<\pi$. 
Let $d'g=\frac{1}{2\pi}|a|^{-2}
d^{*}cd^{*}adxd\theta$ be the measure on $GL_2^+({\Bbb R})$ and thus on 
$GL_2({\Bbb R})$.
Similarly any $g\in SL_2({\Bbb R})$ can be written uniquely as
 $g=n(x)\mat{a}{}{}{a^{-1}}k(\theta)$
with $ x\in {\Bbb R}, a>0$ and $0\leq \theta<2\pi$. Let $d'g=\frac{1}{2\pi}|a|^{-2}
d^{*}adxd\theta$ be the measure on $SL_2({\Bbb R})$.

We now compare the measures $dg$ and $d'g$ on $GL_2$ and $SL_2$ respectively. When
$v$ is a nonarchimedean place for the rational prime $p$, we have $dg=(1+p^{-1})d'g$ in
both case $GL_2$ and $SL_2$. When $v$ is archimedean, we note the measures $d'$
defined
on $GL_2$ and $SL_2$ induce the same quotient measure on $Z({\Bbb R})\backslash 
GL_2({\Bbb R})\cong Z\cap SL_2({\Bbb R})\backslash 
SL_2({\Bbb R})$, and the measures $d$ defined in the introduction also induce the same quotient
measure. Thus our change of the measures is consistent and the equation~(\ref{E:petersson})
still holds.
\end{proof}

\subsection{Ramanujuan conjecture}\label{S:rc}

We show here that the conjecture (\ref{E:rc2}) implies the 
conjecture (\ref{E:rc1}). 

Let $g(z)$ be as in (\ref{E:rc1}). Then $\pp=t(g)$ is a linear combination
of vectors in the space $\tilde{A}_{00}\cap \tilde{A}'(4N,\chi)$ for
some $N$ and $\chi$. We may as well assume
$g(z)$ correspond to a vector $\pp=t(g)$ in a subrepresentation $\ppi$ of
$\tilde{A}_{00}$. From (\ref{E:fc2}), $c(n)=e^{2\pi n} \W^n_{\pp}(e)$.
From the definition,
$$|d_{\ppi}(\pp,S_{\infty},\psi^n)|=
\frac{|\W_{\pp}^n(e)|}{\|\pp\|}e(\pp_{\infty},\psi^n)^{1/2}.$$
From Proposition~\ref{P:86}
$$|d_{\ppi}(\pp,S_{\infty},\psi^n)|=e^{-2\pi n}|c(n)|
[\frac{1}{2} e^{4\pi n}(4\pi n)^{-(1/2+k)}\Gamma(1/2+k)]^{1/2}/\|\pp\|.$$
Thus $|d_{\ppi}(\pp,S_{\infty},\psi^n)|=\delta(\pp)|c(n)|n^{-1/4-k/2}$
where $\delta(\pp)$ is a positive constant depending only on $\pp$. From
(\ref{E:rc2}), 
$$|d_{\ppi}(\pp,S_{\infty},\psi^n)|<<_{\pp,\alpha} n^{\alpha-1/2}.$$
Thus we get $|c(n)|<<_{\pp,\alpha}n^{k/2-1/4+\alpha}$ and (\ref{E:rc1}).

\subsection{Choice of $\pp_2$}

Let $\ppi$ be a cuspidal representation such that the space of $\ppi$ has nontrivial
intersection with $\tilde{A}'_{k+1/2}(4N,\chi)$. 
 The condition (\ref{S:half}.4)  puts a restriction
on $\ppi_2$ (the component at place $v=2$ of $\ppi$). The
 representation $\ppi_2$ must be a subrepresentation of 
$\ppi(\mu\chi_{-1}^{k'},\psi)$ where
$\mu$ is an unramified character of ${\Bbb Z}_2^*$, and $k'=k$ if
$\chi_{(2)}(-1)=1$ and $k'=k+1$ if $\chi_{(2)}(-1)=-1$. 
The space of  vectors in $\ppi(\mu\chi_{-1}^{k'},\psi)$ satisfying (\ref{S:half}.4) is then
two dimensional. It is spanned by two vectors $F[2,1]$ and
$F[2,2^2]$ (\cite{wa1} Proposition~12). 
We recall their definitions (\cite{wa1} p.415, 427). They
are the unique functions in the space of 
$\ppi(\mu\chi_{-1}^{k'},\psi)$ satisfying:
$$F[2,1](\w)=1, \,\,F[2,1](\mat{1}{}{c}{1},1)=0,\,c\in 2{\Bbb Z}_2.
$$
$$F[2,2^2](\w)=0,\,\,F[2,2^2](\mat{1}{}{c}{1},1)=
\text{char}({\Bbb Z}_2)(2^{-2}c).$$

We make a choice of a vector $\pp_2$ in the above two dimensional space.
Define the linear combination
\begin{equation}\label{E:Kvector}
\pp_2=\mu(2^2)\frac{1+(-1)^{k'}i}{4}F[2,1]+F[2,2^2].
\end{equation}
The reason for this choice is explained by the following Proposition.
 Recall the definition of the Whittaker functional $\Lv^z$ by
equation (\ref{E:84}). 
\begin{proposition}\label{C:Kvector}
When
$(-1)^{k'}z\equiv 2,3$ mod $4$, $\L_2^z(\pp_2)=0$.
\end{proposition}
This is a direct consequence of the following computation of
$\L_2^z(F[2,1])$ and $\L_2^z(F[2,2^2])$.
\begin{lemma}\label{L:91}
 With above definitions, $\L_2^z(F[2,1])=\text{char}
({\Bbb Z}_2)(z)$ and
$$\L_2^z(F[2,2^2])=\left\{\begin{array}{ll}
0 &|z|_2>1,\\
(\mu(2^2)+\sqrt{2}
\mu(2^3))\frac{1+(-1)^{k'}i}{4}&z\in (-1)^{k'}+P^2,\\
-\mu(2^2)\frac{1+(-1)^{k'}i}{4}&z\in ((-1)^{k'+1}+P^2)\cup (2+P^2),\\
(\mu(2^2)-\mu(2^4))\frac{1+(-1)^{k'}i}{4}&(-1)^{k'}\frac{z}{4}\in (2+P^2)\cup(-1+P^2).
\end{array}
\right.$$
\end{lemma}
\begin{proof}
The claim for $F[2,1]$ is easy to verify. For $F[2,2^2]$,
using the Iwasawa decomposition, we see:
$$\L_2^z(F[2,2^2])=\int_{|x|_2\geq 2^2}\mu(x^{-1})|x|_2^{-1}
\tilde{\gamma}(x)\chi_{-1}(x)^{k'+1}e^{2\pi izx}dx$$
Consider the integral
$$T(z,i)=\int_{|x|_2= 2^i}\mu(x^{-1})|x|_2^{-1}
\tilde{\gamma}(x)\chi_{-1}(x)^{k'+1}e^{2\pi izx}dx$$
Then \begin{equation}\label{E:91}
\L_2^z(F[2,2^2])=\sum_{i=2}^{\infty}T(z,i).\end{equation}

If $l=2m$ is even, then a change of variable $x\mapsto x2^{-l}$ gives
$T(z,l)=\mu(2^l)T(2^{-l}z,0)$. Over $|x|_2=1$, we have (\cite{wa1}, p. 382)
$$\tilde{\gamma}(x)=1/2(1-i+(1+i)\chi_{-1}(x))$$
Define $\eta(\nu,t)$ to be the Gauss sum: (\cite{wa1}, p.382)
$$\int_{|u|_2=1} \nu(x)e^{-2\pi i tu}d^{*}u$$
Then
$$T(2^{-l}z,0)=(1-2^{-1})^{-1}[\frac{1-i}{2}\eta(\chi_{-1}^{k'+1}, -2^{-l}z)
+\frac{1+i}{2}\eta(\chi_{-1}^{k'}, -2^{-l}z)].$$
Thus
\begin{equation}\label{E:T2}
T(z,2m)=2\mu(2^{2m})[\frac{1-i}{2}\eta(\chi_{-1}^{k'+1}, -2^{-2m}z)
+\frac{1+i}{2}\eta(\chi_{-1}^{k'}, -2^{-2m}z)].
\end{equation}
If $l=2m+1$ is odd, then using the formula $\tilde{\gamma}(2^{-1}x)=\chi_2(x)
\tilde{\gamma}(x)$ and make a change of variable $x\rightarrow
2^{-1}x$, we get $$T(z,2m+1)=\mu(2)\int_{|x|_2=2^{2m}}\mu(
x^{-1})|x|_2^{-1}
\tilde{\gamma}(x)\chi_2(x)\chi_{-1}(x)^{k'+1}e^{\pi izx}dx$$
which by above argument becomes:
\begin{equation}\label{E:T1}
T(z,2m+1)=2\mu(2^{2m+1})[\frac{1-i}{2}\eta(\chi_{-1}^{k'+1}\chi_2, -2^{-2m-1}z)
+\frac{1+i}{2}\eta(\chi_{-1}^{k'}\chi_2, -2^{-2m-1}z)].
\end{equation}
Note
that Gauss sum $\eta(\nu,t)$ vanish if the conductor of $\nu$ is
nonzero and not equal to $-v(t)$, or if $\nu$ is unramified
and $|t|_2>2$. Observe that $\chi_{-1}$ is of conductor 2, and
$\chi_2$ is of conductor 3.
Thus 
$$\L_2^z(F[2,2^2])=\left\{ \begin{array}{ll}
0 &|z|_2>1,\\
T(z,2)+T(z,3)&|z|_2=1,\\
T(z,2)&|z|_2=2^{-1},\\
T(z,2)+T(z,4)+T(z,5)&|z|_2=2^{-2},\\
T(z,2)+T(z,4)&|z|_2=2^{-3}.
\end{array}
\right.$$

We can use the following formulas for $\eta$ (\cite{wa1}, p.383) to finish the
computation: $\eta(\chi_2, 2^{-3})=\frac{1}{\sqrt{2}}$,
$\eta(\chi_{-2}, 2^{-3})=\frac{-i}{\sqrt{2}}$, and
$\eta(\chi_{-1}, 2^{-2})=-i$ (there is a typo in \cite{wa1} for this
value). Note also that $\chi_{-1}(\pm 1+P^2)=\pm 1$ and 
$\eta(\nu,t t')=\eta(\nu,t)\nu^{-1}(t')$ when
$|t'|_2=1$.
Our assertion follows the formulas (\ref{E:T2}) and (\ref{E:T1}). 
\end{proof}

\subsection{The Kohnen space}

Kohnen introduced a subspace $S^+_{k+1/2}(4N,\chi)$ in $S'_{k+1/2}(4N,\chi)$
in \cite{k2}, (we note the notation in \cite{k2} is different from ours).
It consists of $g(z)=\sum_{n=1}^{\infty}c(n)e^{2\pi i nz}$ with
the Fourier coefficients $c(n)$ satisfying:
\begin{equation}\label{E:Kvanish}
c(n)=0, \,\,\text{  when  }\chi_{(2)}(-1)(-1)^kn\equiv 2,3\text{  mod  }4.
\end{equation}

With our definition of $\pp_2$, the Kohnen space has a natural interpretation in the representation
language. Let $\tilde{A}^+_{k+1/2}(4N,\chi)$ 
 be the space generated by vectors 
$\pp=\otimes_v\pp_v$ in the space of cuspidal automorphic forms on $\tilde{SL}_2(\aa)$
satisfying (\ref{S:half}.1)--(\ref{S:half}.4) and with $\pp_2$ being the vector defined in (\ref{E:Kvector}).

\begin{corollary} \label{C:Kspace}
The bijection $g(z)\mapsto t(g)=\pp$ restricts to a bijection between 
 the Kohnen space $S^+_{k+1/2}(4N,\chi)$ and $\tilde{A}^+_{k+1/2}(4N,\chi)$.
\end{corollary}
\begin{proof}
Assume $\pp=t(g)$ with $g(z)=\sum_{n=1}^{\infty}c(n)e^{2\pi inz}$.
From (\ref{E:fc2}), $c(n)=0$ if and only if $\W^n_{\pp}(e)=0$. 

 Let $\pp$ be a vector
as above with $\pp_2$ satisfying (\ref{E:Kvector}). 
Let $n$ be a positive integer such that $\chi_{(2)}(-1)(-1)^{k}n\equiv 2,3$
mod $4$. Since $(-1)^{k'}=\chi_{(2)}(-1)(-1)^k$,
 at $v=2$, $\L_2^n(\pp_2)=0$. From the uniqueness
of the local Whittaker functionals, $\W^n_{\pp}(e)$ vanishes when
$\Lv^n(\pp_v)$ vanish for any place $v$. We get $\W^n_{\pp}(e)=0$ for such
$n$. Thus
 $g(z)=t^{-1}(\pp)$
lies in the Kohnen space, and 
 $\tilde{A}^+_{k+1/2}(4N,\chi)\subset t(S^+_{k+1/2}(4N,\chi))$.

In Proposition~1 of \cite{k2}, Kohnen defined an operator $Q$ on $S'_{k+1/2}(4N,\chi)$.
The operator has two different eigenvalues on this space and $S^+_{k+1/2}(4N,\chi)$ is
the eigenspace of one eigenvalue (denoted $\alpha$). The operator $Q$ induces an operator $Q'$ on the spaces
$V_{\ppi}\cap \tilde{A}'_{k+1/2}(4N,\chi)$. We have a factorization
$V_{\ppi}\cap \tilde{A}'_{k+1/2}(4N,\chi)=\otimes V'_{\ppi,v}$ with $V'_{\ppi,2}$ a
two dimensional space.  Then $Q'=\otimes Q'_v$. In fact $Q'_v$
are all trivial actions for $v\not=2$.   Clearly $\pp_2$ in (\ref{E:Kvector}) 
is the eigenvector of $Q'_2$ with
eigenvalue $\alpha$ as the  vector $\pp=\otimes\pp_v$ with local component $\pp_2$ lies in 
$\tilde{A}^+_{k+1/2}(4N,\chi)$. 
Let $\tilde{A}^-_{k+1/2}(4N,\chi)$ be the subspace of $\tilde{A}'_{k+1/2}(4N,\chi)$
generated by $\pp'=\otimes\pp'_v$ where $\pp_2'$ is 
the eigenvector for the other eigenvalue. Then $\tilde{A}'_{k+1/2}(4N,\chi)=
\tilde{A}^+_{k+1/2}(4N,\chi)\oplus
\tilde{A}^-_{k+1/2}(4N,\chi)$. 
As $\tilde{A}^-_{k+1/2}(4N,\chi)\cap  t(S^+_{k+1/2}(4N,\chi))=0$,
we get the corollary.
\end{proof}

\subsection{Local computation at the place $2$}

We compute the local factor as in \S~\ref{S:local}. The vector $\varphi_2$ is the unramified
vector chosen as in subsection~\ref{S:principal}.
\begin{proposition}\label{P:e2}
 Assume
$\mu_2(x)=|x|_2^{ir}$, with $r\in {\Bbb R}$. Then
$$e(\varphi_2,\psi)=3/2 |1-q^{-1-2ir}|^{-2},$$
 $$e(\pp_2,\psi^{|D|})=\left\{\begin{array}{ll}
3/4|1+2^{-1/2-ir}|^{-2}& D\in 1+P^2,\\
3/4|1-2^{-1-2ir}|^{-2}|D|_2^{-1}& \frac{D}{4}\in (2+P^2)\cup (-1+P^2).
\end{array}\right.$$

Therefore, when $D\in 1+P^2$,
$$\frac{e(\varphi_2,\psi)}{e(\pp_2,\psi^{|D|})}=\left\{\begin{array}{ll}
2|D|_2 L(\pi_2,1/2)& D\in 1+P^2,\\
2|D|_2& \frac{D}{4}\in (2+P^2)\cup (-1+P^2).
\end{array}\right.$$
\end{proposition}
\begin{proof}
The formula for $e(\varphi_2,\psi)$ is given in Proposition~\ref{P:81}.
One can use (\ref{E:norm6}) to compute $\|\pp_2\|$ when $\mu_2(x)=|x|^{ir}$
where $r\in {\Bbb R}$.
Since $F[2,1]$ and $F[2,2^2]$ are perpendicular, we get
$$\|\pp_2\|^2=1/8\|F[2,1]\|^2+\|F[2,2^2]\|^2$$
Use the Iwasawa decomposition it is easy to show 
$$\|\pp_2\|^2=(1/4+1/8)|D|_2^{-1}=3/8|D|_2^{-1}.$$
From Lemma~\ref{L:91}, we get the claim for $e(\pp_2, \psi^{|D|})$.
From the formula on local $L-$factor in \cite{g}, we get the last statement of the Proposition.
\end{proof}

\section{A generalization of the Kohnen-Zagier formula}\label{S:KZ}

\subsection{Statement of the Theorem}

Let $f(z)$ be a cusp form as in (\ref{E:11}) with square free level
$N$ (odd) and weight $2k$. Let $S_N$ be the set of primes $p|N$. 
Let $S$ be a (possibly empty)  subset of $S_0$. Define $D_S$ to be the
set of  fundamental
discriminants $D$ such that $(\frac{D}{p})=-w_p$ if and only if $p\in S$. Then the
set of fundamental discriminants is the disjoint union $\cup_{S\subset S_N}D_S$.
For $D$ a fundamental discriminant, let
 $T(D)$ be the set of $p|N$ that also divides $D$. Let 
$sgn(D)=D/|D|$. The character $\psi$ is defined as in \S~\ref{S:dictionary}.

\begin{theorem}\label{T:KZgeneral}
Let $S\subset S_N$ and $s$ be the size of $S$.
Let $N'=\prod_{p\in S}p$,
let $\chi=\prod_{p|2N}\chi_{(p)}$ be any Dirichlet character of $({\Bbb Z}/4NN')^*$ such that
$\chi_{(p)}\equiv 1$ when $pN'|N$, $\chi_{(p)}(-1)=-1$ when $p| N'$ and $\chi(-1)=1$. 
There exists a unique (up to scalar multiple)
 cusp form $g_S(z)$ in $S'_{k+1/2}(4NN',\chi)$,
such that the following is true:

(1) $g_S(z)$ is a Shimura lift of $f(z)$.

(2) $g_S(z)$ lies in the Kohnen space, i.e. if $g_S(z)=\sum_{n=1}^{\infty}c(n)e^{2\pi inz}$,
then $c(n)=0$ when $(-1)^{s+k}n\equiv 2,3$ mod 4.

(3) $c(D)= 0$ if $D>0$ is a fundamental discriminant with $(-1)^{s+k}D\not\in D_S$.

Moreover for this $g_S(z)$ and for $D\in D_S$, 
if $(-1)^{s+k}\not=sgn(D)$, then 
$L(f,D,k)=0$; 
if $(-1)^{s+k}=sgn(D)$, then 
\begin{equation}\label{E:kzg}
\frac{|c(|D|)|^2}{<g_S,g_S>}=\frac{L(f,D,k)}{<f,f>}
|D|^{k-1/2}\frac{(k-1)!}{\pi^k}2^{\nu(N)-t}
\prod_{p\in S}\frac{p}{p+1}.
\end{equation}
Here   $t$ is the size of $T(D)$.
\end{theorem}

From subsection~\ref{S:integer}, the new form $f(z)$ determines an irreducible cuspidal
representation $\pi$ of $GL_2(\aa_{\Bbb Q})$ with trivial central character. 
Here we say $g(z)$ is a Shimura lift of $f(z)$ if $\pp=t(g)$ lies in a space $V_{\ppi}$
where $\ppi$ satisfies $\pi=S_{\psi}(\ppi)$ or $\pi=S_{\psi^{-1}}(\ppi)$
(see Theorem~\ref{T:Wglobal} for the notion $S_{\psi}$).

The proof of the Theorem involves  translating Theorem~\ref{T:KZadele} into the language
of cusp forms. We will set up the translation in subsections (\ref{S:KZ2})--(\ref{S:KZ3}).

\subsection{A Lemma on Atkin-Lehner involution}\label{S:KZ2}

 Let $\hat{w}_p=\mat{p}{a}{N}{pb}$ with $a,b \in {\Bbb
Z}$ and $\det \hat{w}_p=p$. Let
$\tilde{f}(z)=f|_{p^{-1/2}\hat{w}_p}(z)$ be the Atkin-Lehner
involution. Let $\hat{\varphi}=s(\tilde{f})$. For lack of
reference, we give a proof of the following well known result.
\begin{lemma}\label{L:94}
With the above definition, then $\pi_p=\sigma(\mu,\mu^{-1})$
with $\mu(x)=|x|^{1/2}\chi_{\tau}(x)$, $\tau$ a unit
in ${\Bbb Z}_p$. We have $\tilde{f}=w_p f$, with
$w_p=1$ when $\tau$ is not a square, and $w_p=-1$ when
$\tau$ is a square. Moreover $\epsilon(\pi_p,1/2)=w_p$
\end{lemma}
\begin{proof}
 The first claim is in \cite{ge}. Since
$\hat{\varphi}(g_{\infty})=\varphi(\hat{w}_p g_{\infty})$, from
left $GL_2({\Bbb Q})Z(\aa_{\Bbb Q})$ invariance,
$$\hat{\varphi}(g_{\infty})=\varphi(g_{\infty}\prod_{v\not=\infty}
\hat{w}_{p,v}^{-1})=\pi(\prod_{v\not=\infty}
\hat{w}_{p,v}^{-1})\varphi(g_{\infty})$$
As $\varphi=\otimes \varphi_v$, we get $\hat{\varphi}=
\varphi_{\infty}\otimes_{v\not=\infty}\pi_v(\hat{w}_{p,v}^{-1})
\varphi_v$. When $v$ does not divide $N$, $\hat{w}_{p,v}^{-1}
\in GL_2({\Bbb Z}_v)$, thus fixes $\varphi_v$. When $v|N$
but $v\not=p$, $\hat{w}_{p,v}^{-1}\in K_{0,v}$, thus fixes
$\varphi_v$. When $v=p$, $\hat{w}_{p,p}(w\u{p})^{-1}\in K_{0,p}$,
thus $\pi_p(\hat{w}_{p,p}^{-1})\varphi_{p}=\pi_p(w\u{p})^{-1}
\varphi_{p}$. As $w\u{p}K_{0,p}(w\u{p})^{-1}\in K_{0,p}$,
the vector $\pi_p(w\u{p})^{-1}
\varphi_{p}$ is again
fixed by $K_{0,p}$, thus is a scalar multiple of $\varphi_{p}$.
Denote this multiple by $w_p$. Then $\hat{\varphi}=
w_p\varphi$, thus $\tilde{f}=w_p f$.
To find the multiple, we only need to evaluate
$\pi_p(w\u{p})^{-1}
\varphi_{p}(e)$ which clearly equals $p\chi_{\tau}(p)$.
Since $\varphi_{p}(e)=-p$, we get $w_p=-\chi_{\tau}(p)$
which gives the claim in the Lemma. The computation of
$\epsilon(\pi_p,1/2)$ is given in \cite{g}.
\end{proof}

\subsection{Definition of $\epsilon(S)$}

 Let $f$ and $\pi$ be
as before.
Since $\pi_v$ is unramified for all places $v$ where
$v\not=\infty$ and
$|N|_v=1$, we can let $\Sigma$ in Theorem~\ref{T:Wglobal}
to be the set $\{\infty\}\cup S_N$.
Let $S$ be a set as in the Theorem. Then it determines an $\epsilon(S)\in
\{\pm 1\}^{|\Sigma|}$, where the component  $\epsilon(S)_p$
at $p\in S$ is $-w_p$, at $p|N$ and $p\not\in S$ is $w_p$, and
 at $\infty$ is $(-1)^{s+k}$. 

\begin{lemma}\label{L:p1}
We have $\epsilon(\pi,1/2)=\prod_{v\in \Sigma}\epsilon(S)_v$.
\end{lemma}
\begin{proof}
The product on the right is $(-1)^k\prod_{p|N}w_p$. Since $w_p=\epsilon(\pi_p,1/2)$
by Lemma~\ref{L:94} and $(-1)^k=\epsilon(\pi_{\infty},1/2)$ from
\cite{g}, we get the claim.
\end{proof}

From Theorem~\ref{T:Wglobal}, associated to $\pi$ and the character
$\psi_S(x)=\psi((-1)^{k+s}x)$ is the Shimura lift
$\ppi^{\epsilon(S)}$. Here $\ppi^{\epsilon(S)}=\Theta(\pi\otimes\chi_D,\psi_S^D)$
for some $D\in {\Bbb Q}^{\epsilon(S)}(\pi)$.

\subsection{Relation between $D_S$ and ${\Bbb Q}^{\epsilon(S)}(\pi)$}

Given $D$ a fundamental discriminant, let 
$\epsilon_v(D)=(\frac{D}{\pi_v})$ for $v\in \Sigma$.

\begin{lemma}\label{L:epsilon}
When $v=\infty$, $\pi_{\infty}$ is a discrete series,
$\epsilon_{\infty}(D)=sgn(D)$. 

When $p|N$, $\pi_p$ is a special
representation $\sigma^{\tau_v}$ as in subsection~\ref{S:special}, where
$\tau_v$ is a unit in ${\Bbb Z}_p$.  Then  $\epsilon_p(D)=w_p$ if $p|D$;
$\epsilon_p(D)=(\frac{D}{p})$ when $D$ is a unit in  ${\Bbb Z}_p$.
\end{lemma}
\begin{proof}
At $v=\infty$, $\pi_{\infty}\cong \pi_{\infty}\chi_D$,  thus $\epsilon_{\infty}(D)
=\chi_{D}(-1)=sgn(D)$. 
When
$p|N$, as $\Theta(\pi_p,\psi)$ is either a special representation (when
$\tau_v$ is not a square) or an odd Weil representation (when $\tau_v$ is
a square). 
Note that  $w_p=1$ if and only if $\tau$ is not a square 
(see Lemma~\ref{L:94}).
In the case $\tau_v$ is not a square, $\Theta(\pi_p,\psi)$ has
a nontrivial $\psi^D-$Whittaker model when $D$ is not a nonsquare unit.
In the case
$\tau_v$ is a square, then only when $D$ is a square does $\Theta(\pi_p,\psi)$
has a nontrivial $\psi^D-$Whittaker model.
From
Theorem~\ref{T:Wlocal},  we get the result. 
\end{proof}

\begin{lemma}\label{L:p2}
When $D\in D_S$, $\epsilon_p(D)=\epsilon(S)_p$ for all $p|N$. 
\end{lemma}
\begin{proof}
When $p\in S$, then $D$ is a unit in ${\Bbb Z}_p$, and
 $\epsilon_p(D)=(\frac{D}{p})=-w_p=\epsilon(S)_p$.
When $p\in S_N-S$, then either $p|D$ in which case $\epsilon_p(D)=w_p=\epsilon(S)_p$
or $D$ is a unit in ${\Bbb Z}_p$, in which case
$\epsilon_p(D)=(\frac{D}{p})=w_p=\epsilon(S)_p$.
\end{proof}

As the set ${\Bbb Q}^{\epsilon(S)}(\pi)$ consists of $D$ with 
$\epsilon_v(D)=\epsilon(S)_v$ for $v\in \Sigma$, from the above lemma and the formula
for $\epsilon_{\infty}(D)$, we get
\begin{corollary}
A fundamental discriminant $D$ lies in ${\Bbb Q}^{\epsilon(S)}(\pi)$ if
and only if $D\in D_S$ and $(-1)^{s+k}=sgn(D)$.
\end{corollary}

\subsection{Description of $g_S$}\label{S:KZ3}

The cusp forms $g_S$ in the Theorem is taken to be the inverse image
$t^{-1}(\pp_S)$ of some vector $\pp_S$ in the space of $\ppi^{\epsilon(S)}$.
We describe the choice of $\pp_S=\otimes\pp_v$.

Using the explicit description of theta correspondence in \cite{wa2}, we
get the following description on the local components of $\ppi^{\epsilon(S)}
=\otimes \ppi^{\epsilon(S)}_v$. Note that $\ppi^{\epsilon(S)}_v
=\Theta(\pi_{v}\otimes\chi_{D},\psi_S^{D})$ for some $D\in {\Bbb Q}^{\epsilon(S)}(\pi)$.

Recall that the description of $\pi=\otimes\pi_v$ 
is given in (\ref{S:integer}.1)--(\ref{S:integer}.3),
along with a choice of the vector $\varphi=\otimes\varphi_v$
such that $\varphi=s(f)$. Below is the description of $\pi_v$, $\ppi^{\epsilon(S)}_v$ and
the choice of the vectors $\pp_v$.

(\ref{S:KZ3}.1) When
$v=\infty$, $\pi_v$ is the discrete series
$\sigma(\mu_{\infty},\mu_{\infty}^{-1})$, 
with $\mu_{\infty}(x)=|x|_v^{k-1/2}(\text{sgn }x)^k$.
When $\epsilon_{\infty}(D)=\epsilon(S)_{\infty}$,
we get $sgn(D)=(-1)^{s+k}$ thus
$\psi_S^D=\psi^{|D|}$. 
Thus $\ppi^{\epsilon(S)}_{\infty}=
\Theta(\pi_{\infty}\otimes\chi_{D},\psi^{|D|})$.
As $\pi_{\infty}\otimes\chi_{|D|}\cong \pi_{\infty}\otimes \chi_D$ and 
$|D|$ and $1$ are in the same square class, we get $\ppi^{\epsilon(S)}_{\infty}=
\Theta(\pi_{\infty},\psi)=\tilde{\sigma}_{\infty}(\mu_{\infty})$, (\cite{wa2}).
 We take
$\pp_{\infty}$ to be the vector with minimal weight.

(\ref{S:KZ3}.2) At $p\not\in \Sigma$, 
$\pi_p=\pi(\mu_p, \mu_p^{-1})$ with $\mu_p$ an unramified
character. From a well known result of Deligne,  the
unramified characters $\mu_p$ has the form $\mu_p(x)=|s|^{ir}$ with $r\in {\Bbb R}$,
(this is the Ramanujuan conjecture for the integral weight forms).
 Then
$\ppi^{\epsilon(S)}_p=\Theta(\pi_p,\psi_S)=\ppi(\mu_p\chi_{-1}^{s+k},\psi)$. 
We take $\pp_p$
to be the unramified vector in this unramified representation when $p\not=2$. 
We let 
$\pp_2$ be the vector defined by (\ref{E:Kvector}) with $k'=k+s$.

(\ref{S:KZ3}.3) At $p\in S_N-S$,  $\pi_p=\sigma^{\tau}$ 
with $\tau\in {\Bbb Z}_p^*$. Let $D$ be a unit in ${\Bbb Z}_p$
such that $\tau D$ is not a square,
then $\epsilon_p(D)=w_p=\epsilon(S)_p$.
Thus $\ppi^{\epsilon(S)}_p=\Theta(\sigma^{\tau D},\psi_S^{D})= 
\tilde{\sigma}^{\delta }(\psi_S^{D})$; here $\delta$ is
any nonsquare unit.
  We take $\pp_p$ to be the vector in Lemma~\ref{L:84}.

(\ref{S:KZ3}.4) When  $p\in S$, again $\pi_p=\sigma^{\tau}$ 
with $\tau\in {\Bbb Z}_p^*$. Let $D$ be a unit in ${\Bbb Z}_p$
such that $\tau D$ is  a square, then $\epsilon_p(D)=-w_p=\epsilon(S)_p$.
Thus $\ppi^{\epsilon(S)}_p=\Theta(\sigma^{1},\psi_S^{D})$ which
is 
$r^-_{\psi_S^{D}}$ from subsection~\ref{S:special}. Let
$\chi_{(p)}$ be the character on ${\Bbb Z}_p^*$ defined in the Theorem, we let
$\pp_p=\Phi_{\chi_{(p)}}$ where $\Phi_{\chi_{(p)}}$ is defined in Lemma~\ref{L:weil1}.

Each of the choices of $\pp_p$ is determined unique up to a scalar
multiple. Let $\pp=\otimes \pp_v$. We define the cusp form $g_S(z)$ to be $t^{-1}(\pp)$.

\subsection{Proof of the Theorem}

\begin{proof}
As $S_{\psi_S}(\ppi^{\epsilon(S)})=\pi$, we see
 $g_S(z)$ is a Shimura lift of $f$.
We can check $\pp=\otimes\pp_v\in \tilde{A}'_{k+1/2}(4NN',\chi)$.
Thus $g_S(z)\in S'_{k+1/2}(4NN',\chi)$.
It lies in the Kohnen space
because of our choice of $\pp_2$. If $D$ is a fundamental discriminant with $\pm D\not\in D_S$,
then $\pm D\not\in {\Bbb Q}^{\epsilon(S)}(\pi)$. By Theorem~\ref{T:Wadele}, we get
$d_{\ppi^{\epsilon(S)}}
(\Sigma\cup\{2\},\psi_S^{\pm D})=0$. As $\psi^{|D|}$ is one of $\psi_S^{\pm D}$,
we get $c(|D|)= 0$ from the consideration in \S~\ref{S:dictionary}.

Next we show the uniqueness. If $g(z)\not=0$ is a Shimura lift of $f$ such that $t(g)$ lies in the space
of $\ppi$, we show $\ppi=\ppi^{\epsilon(S)}$. As $\ppi_{\infty}$ is a holomorphic discrete 
series, by Theorem~\ref{T:Wglobal} $\ppi=
\Theta(\pi\otimes\chi_{D_1},\psi^{|D_1|})$  for some $D_1$.
As $g(z)\not=0$, $c((-1)^{s+k}D_2)\not=0$ for some $D_2\in D_S$. The condition 
  implies that $D_2\in {\Bbb Q}^{\epsilon(S)}(\pi)$ and $\ppi_v$
has nontrivial $\psi^{|D_2|}$-Whittaker model at all places $v$. 
From Theorem~\ref{T:Wlocal},
we see $\ppi=\Theta(\pi\otimes\chi_{\alpha D_2}, \psi^{|D_2|})$, 
where $\alpha=\frac{|D_1D_2|}{D_1D_2}=\pm 1$.
 Examine the component $\ppi_2$. As $\chi_{(2)}(-1)=(-1)^s$ under
our assumptions, we get $\ppi_2=\Theta(\pi\otimes\chi_{-1}^{s+k},\psi)$, thus
we see $\ppi=\Theta(\pi\otimes\chi_{D_2},\psi^{|D_2|})=\ppi^{\epsilon(S)}$. From
Corollary~\ref{C:Kspace} and the fact that $\tilde{A}^+_{k+1/2}(4N,\chi)\cap 
V_{\ppi^{\epsilon(S)}}$ is one dimensional, we get the uniqueness of $g_S(z)$.

We now prove the identity (\ref{E:kzg}).
Let $D\in D_S$. 
If $(-1)^{s+k}\not=sgn(D)$, then from equation (\ref{E:p3}),
Lemma~\ref{L:p1}
and \ref{L:p2}, we get:
$$\epsilon(\pi\otimes\chi_D,1/2)=\epsilon(\pi,1/2)\prod_{v\in\Sigma}\epsilon_v(D)=
-\epsilon(\pi,1/2)\prod_{v\in\Sigma}\epsilon(S)_v=-1.$$
Thus $L(\pi\otimes\chi_D,1/2)=0$, i.e. $L(f,D,k)=0$.

When $(-1)^{s+k}=sgn(D)$, then 
$D\in {\Bbb Q}^{\epsilon(S)}(\pi)$, thus we can apply
Theorem~\ref{T:KZadele} to get:
\begin{equation}\label{E:kz2}
L^{\Sigma\cup\{2\}}(\pi\otimes\chi_D,1/2)=\frac{|d_{\ppi^{\epsilon(S)}}
(\Sigma\cup\{2\},\psi_S^{D})|^2}{|d_{\pi}(\Sigma\cup\{2\},\psi_S)|^2}
\prod_{v\in\Sigma\cup\{2\}}|D|_v.
\end{equation}
From Lemma~\ref{L:dpi}, $d_{\pi}(\Sigma\cup\{2\},\psi_S)=d_{\pi}
(\Sigma\cup\{2\},\psi)$. Observe also $\psi_S^D=\psi^{|D|}$.
From the explicit description of $d_{\ppi^{\epsilon(S)}}
(\Sigma\cup\{2\},\psi^{|D|}) $ and $d_{\pi}
(\Sigma\cup\{2\},\psi)$, we get:
\begin{equation}\label{E:kz3}
\frac{|d_{\ppi^{\epsilon(S)}}
(\Sigma\cup\{2\},\psi^{|D|})|^2}{|d_{\pi}(\Sigma\cup\{2\},\psi)|^2}=
\frac{|\W^{|D|}_{\pp}(e)|^2\|\varphi\|^2}
{|W_{\varphi}(e)|^2\|\pp_S\|^2}
\prod_{v\in \Sigma\cup\{2\}}\frac{e(\pp_v,\psi^{|D|})}{e(\varphi_v,\psi)}.
\end{equation}
Recall the results on the local factors from Propositions~\ref{P:86}, \ref{P:85},
 \ref{P:weil1} and \ref{P:e2}:
\begin{equation}\label{E:ep}
\frac{e(\varphi_p,\psi)}
{e(\pp_{p},\psi^{|D|})}=\left\{\begin{array}{ll}
\frac{1}{2}e^{4\pi (1-|D|)}|D|^{1/2+k}\pi^{-k}(k-1)! & p=\infty,\\
2L(\pi_p\otimes\chi_D,1/2)|D|_p&p\in S_N-S, p\not| N,\\
L(\pi_p\otimes\chi_D,1/2)|D|_p&p\in S_N-S, p|N,\\
2(1+q^{-1})^{-1}L(\pi_v\otimes\chi_D,1/2)&p\in S,\\
2|D|_2L(\pi_2,1/2)& p=2, D\in 1+P^2,\\
2|D|_2& p=2, \frac{D}{4}\in (2+P^2)\cup (-1+P^2).
\end{array}\right.\end{equation}

From subsection~\ref{S:integer}, $W_{\varphi}(e)=e^{-2\pi}$, and from
(\ref{E:fc2}), $\W^{|D|}_{\pp}(e)=e^{-2\pi |D|}c(|D|)$. Thus we get
\begin{equation}\label{E:kz4}
L^{{\infty}\cup\{2\}}(\pi\otimes\chi_{D},1/2)l_2
=(1/2)^{\nu(N)-t}|D|^{-k+1/2}\frac{\pi^k}{(k-1)!}\frac{|c(|D|)|^2\|\varphi\|^2}
{\|\pp_S\|^2}
\prod_{p\in S}(1+p^{-1}).
\end{equation}
Here we set $l_2$ to be  $L(\pi_2,1/2)$ 
when $D\equiv 1$
mod $4$ and to be  $1$ when $D\equiv 0$ mod $4$.
\begin{remark}
We now notice some differences between $L(\pi\otimes\chi_D,s)$
and $L(f,D,s')$. First $L(f,D,s')$ does not have factor at $\infty$.
Secondly, because $\chi_D(2)$ over $v=2$ 
is not the same as $(\frac{D}{2})$, we need to correct the factor at
the place $v=2$. We have actually
$$L(f,D,k)=L^{{\infty}\cup{2}}(\pi\otimes\chi_D,1/2)l_2.$$ 
\end{remark}

From  the above remark, Lemma~\ref{L:petersson} and (\ref{E:kz4}), we get:
\begin{equation}\label{E:kz6}
L(f,D,k)
=(1/2)^{\nu(N)-t}|D|^{-k+1/2}\frac{\pi^k}{(k-1)!}\frac{|c(|D|)|^2<f,f>}
{<g_S,g_S>}.
\end{equation}
 Therefore we get (\ref{E:kzg}).
\end{proof}

\subsection{Some examples}

Example 1: When $S$ is empty, $g_S(z)$ is the $g(z)$ in (\ref{E:11}). If $T(D)$ is also
empty, we recover (\ref{E:11}). If $T(D)$ is nonempty, (\ref{E:11}) should be revised to:
$$\frac{|c(|D|)|^2}{<g,g>}=|D|^{k-1/2}\frac{(k-1)!}{\pi^k}2^{\nu(N)-t}\frac{L(f,D,k)}{<f,f>}.
$$

Example 2: Look at the case when $f(z)$ is the new form of weight 2
and level 11. Such a form exists and is unique, with $w_{11}=-1$. The theorem says there is
a cusp form $g_{\{11\}}(z)$ in $S'_{3/2}(484,\chi)$, where $\chi$ satisfies the
condition in the Theorem, (for example $\chi(x)=(-1)^{(x-1)/2}(\frac{x}{11})$),
such that $c(n)=0$ whenever $n\equiv 2,3$ mod $4$,
and when the fundamental discriminant $D>0$ is such that $(\frac{D}{11})=1$,
$$\frac{|c(D)|^2}{<g_{\{11\}},g_{\{11\}}>}=|D|^{1/2}\frac{11}{6\pi}
\frac{L(f,D,k)}{<f,f>}.$$

\end{document}